\documentclass[11pt]{article}
\usepackage[all]{xy}
\usepackage{amsfonts}
\usepackage{blkarray}
\usepackage{mathtools}
\usepackage{color}
\usepackage{tikz-cd}
\usepackage{float}
\usepackage{mathtools}
\usepackage{amsthm}
\usepackage[margin=1.1in]{geometry}

\providecommand{\image}{\mathop{\rm im}\nolimits}
\providecommand{\kernel}{\mathop{\rm ker}\nolimits}
\providecommand{\trace}{\mathop{\rm tr}\nolimits}
\providecommand{\rank}{\mathop{\rm rank}\nolimits}
\providecommand{\lcm}{\mathop{\rm lcm}\nolimits}

\usetikzlibrary{chains}
\tikzset{node distance=2em, ch/.style={circle,draw,on chain,inner sep=2pt},chj/.style={ch,join},every path/.style={shorten >=4pt,shorten <=4pt},line width=1pt,baseline=-1ex}
\newcommand{\mnode}[2][chj]{%
\node[#1,label={below:{#2}}] {};
}
\newcommand{\dydots}{%
\node[chj,draw=none,inner sep=1pt] {\dots};
}
\newcommand{\mnoder}[2][chj]{%
\node[#1,label={right:{#2}}] {};
}
\definecolor{grey}{gray}{0.5}
\newtheorem{theorem}{Theorem}[section]
\newtheorem{lemma}[theorem]{Lemma}
\newtheorem{proposition}[theorem]{Proposition}
\newtheorem{corollary}[theorem]{Corollary}
\newtheorem{definition}[theorem]{Definition}
\newtheorem{conjecture}[theorem]{Conjecture}
\begin{document}
\title{Linear relations and integrability for cluster algebras from affine quivers}
\author{Joe Pallister}
\date{\vspace{-5ex}}
\maketitle
\begin{abstract}
We consider frieze sequences corresponding to sequences of cluster mutations for affine $D$ and $E$ type quivers. We show that the cluster variables satisfy linear recurrences with periodic coefficients, which imply the constant coefficient relations found by Keller and Scherotzke. Viewing the frieze sequence as a discrete dynamical system, we reduce it to a symplectic map on a lower dimensional space and prove Liouville integrability of the latter.
\end{abstract}
\section{Introduction}\label{generalapproach}
A cluster algebra is an algebra with a distinguished set of generators, cluster variables, that appear in sets called seeds, \cite{clusteri}. Each seed is obtained recursively, by a process called mutation, from an initial seed. In our work there is a quiver associated with each seed, with $N+1$ vertices, without loops or $2-$cycles. It is a theorem of \cite{clusterii} that there are finitely many cluster variables if and only if the quiver associated to one of the seeds is of simply laced Dynkin type. 

From a quiver without cycles we construct a frieze \cite{frises} by a sequence of mutations at only sinks or sources. These mutations give a coupled system of $N+1$ recurrences
\begin{equation}\label{generalfrieze}
X^k_{n+1}X^k_n=1+\left(\prod_{i\rightarrow k}X_n^i\right)\left(\prod_{i\leftarrow k}X_{n+1}^i\right)
\end{equation}
which define the frieze sequence. Here the products are taken over neighbours with arrows in to and out of $k$ respectively in the initial quiver. The superscript runs over each $k=1,\ldots,N+1$, and the subscript $n$ denotes the iterates of the recurrence. We consider the frieze variables with superscript $k$ to be ``living" at vertex $k$. To see the frieze as a discrete dynamical system we set initial variables $X_0^k$ for each vertex $k$, and define recursively $X_{n+1}^k$ by Equation \ref{generalfrieze}. The map
\begin{equation}\label{generalrecursion}
\varphi:
\begin{pmatrix}
X^1_n \\
X^2_n \\
\vdots \\
X^{N+1}_n
\end{pmatrix}
\mapsto
\begin{pmatrix}
X^1_{n+1} \\
X^2_{n+1} \\
\vdots \\
X^{N+1}_{n+1}
\end{pmatrix}
\end{equation}
is called the cluster map. The form of cluster mutation means that each new variable can be written as a birational function of the previous ones, and since (\ref{generalfrieze}) is obtained as a composition of mutations the map (\ref{generalrecursion}) is birational. This induces an automorphism $\varphi^*$, a ``shift", on the field of rational functions
\[
\mathbb{C}(\{X^k_0\}_{k=1,\ldots,N+1})
\] 
by $\varphi^*(X^k_n)=X^k_{n+1}$ for each $n\in\mathbb{Z}$. We say an element of this field is period $p$ if it is fixed by $(\varphi^*)^p$. Nonconstant period one elements are called invariants or first integrals. 

\begin{figure}
\begin{align*}
&\tilde{A}_N &&
\begin{tikzpicture}[start chain,node distance=1ex and 2em]
\mnode{1}
\mnode{1}
\dydots
\mnode{1}
\mnode{1}
\begin{scope}[start chain=br going above]
\chainin(chain-3);
\node[ch,join=with chain-1,join=with chain-5,label={above:{1}}] {};
\end{scope}
\end{tikzpicture}
\\
&\tilde{D}_N (N \ge 4) &&
\begin{tikzpicture}
\begin{scope}[start chain]
\mnode{1}
\mnode{2}
\mnode{2}
\dydots
\mnode{2}
\mnode{1}
\end{scope}
\begin{scope}[start chain=br going above]
\chainin(chain-2);
\mnoder{1};
\end{scope}
\begin{scope}[start chain=br2 going above]
\chainin(chain-5);
\mnoder{1};
\end{scope}
\end{tikzpicture}
\\
&\tilde{E}_6 &&
\begin{tikzpicture}
\begin{scope}[start chain]
\foreach \dyi in {1,2,3,2,1} {
	\mnode{\dyi}
}
\end{scope}
\begin{scope}[start chain=br going above]
\chainin(chain-3);
\mnoder{2}
\mnoder{1}
\end{scope}
\end{tikzpicture}
\\
&\tilde{E}_7 &&
\begin{tikzpicture}
\begin{scope}[start chain]
\mnode{1}
\mnode{2}
\mnode{3}
\mnode{4}
\mnode{3}
\mnode{2}
\mnode{1}
\end{scope}
\begin{scope}[start chain=br going above]
\chainin(chain-4);
\mnoder{2}
\end{scope}
\end{tikzpicture}
\\
&\tilde{E}_8 &&
\begin{tikzpicture}
\begin{scope}[start chain]
\mnode{1}
\mnode{2}
\mnode{3}
\mnode{4}
\mnode{5}
\mnode{6}
\mnode{4}
\mnode{2}
\end{scope}
\begin{scope}[start chain=br going above]
\chainin(chain-6);
\mnoder{3}
\end{scope}
\end{tikzpicture}
\end{align*}\caption{The simply laced affine diagrams with each vertex $i$ labelled $\delta_i$.}\label{affinediagrams}
\end{figure}
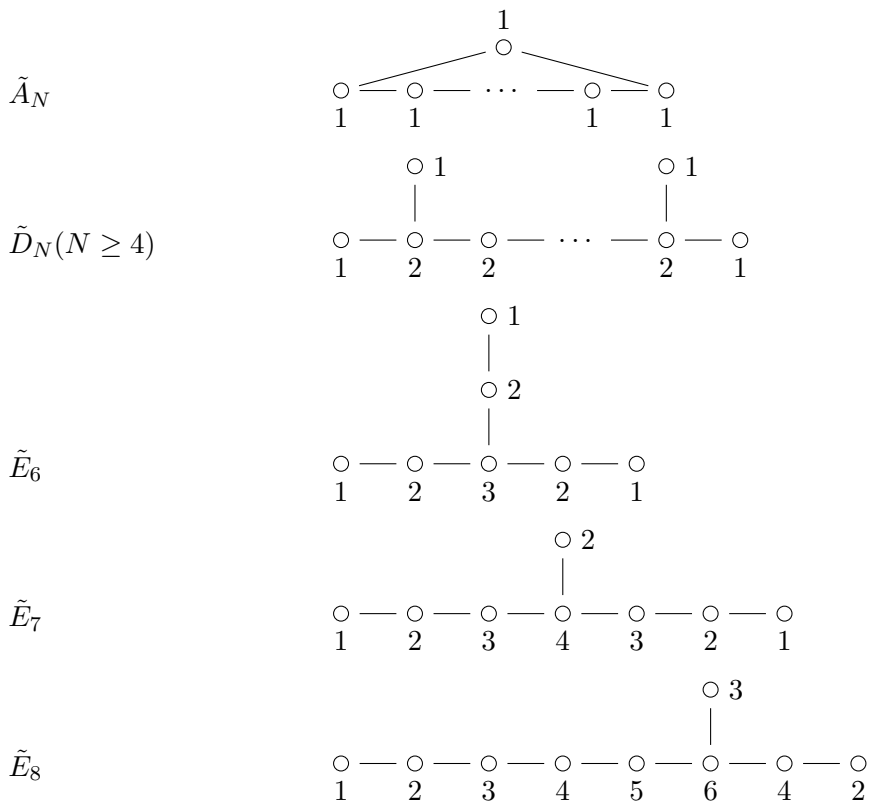

It was proved in \cite{frises} that if there are linear relations between the variables in the frieze sequence then the corresponding diagram is either affine or Dynkin and is simply laced. These are shown in Figure \ref{affinediagrams}
where, and throughout this paper, we have used a tilde to denote the affine version of the diagram or quiver. The converse for $\tilde{A}$ and $\tilde{D}$ types was also proved in \cite{frises}, that the associated frieze sequences satisfy linear relations. In \cite{kellerscherotzke}, using a representation theoretic approach, the authors proved this for all $\tilde{A}\tilde{D}\tilde{E}$ quivers, with some of these linear relations given explicitly. We adopt a different method of proof here, allowing us a uniform expression for these results:
\begin{theorem}
For each extending vertex $k$, the variables that live there satisfy the constant coefficient linear relation
\begin{equation}\label{KSlinearrelation}
X^k_{n+2b}-\mathcal{K}X^k_{n+b}+X^k_n=0.
\end{equation}
Here $b$ depends on the quiver, the values of which are given in Figure \ref{btable}, and $\mathcal{K}$ is invariant under the shift $\varphi^*$. The extending vertices for an affine quiver are those labelled with a $1$ in Figure \ref{affinediagrams}.
\end{theorem}
\begin{figure}
\begin{center}
\begin{tabular}{|r|l|}
\hline 
Quiver & $b$ \\
\hline\noalign{\smallskip}
$\tilde{A}_{p,q}$ & $\lcm(p,q)$\\
$\tilde{D}_N \quad N \text{ even}$ & $N-2$  \\
$\tilde{D}_N \quad N \text{ odd}$ & $2N-4$  \\
$\tilde{E}_6$ & $6$ \\
$\tilde{E}_7$ & $12$ \\
$\tilde{E}_8$ & $30$ \\
  \hline
\end{tabular}
\end{center}\caption{Values of $b$ for the $\tilde{A}\tilde{D}\tilde{E}$ quivers.}\label{btable}
\end{figure}
We remark that in \cite{kellerscherotzke}, for $\tilde{D}$ and $\tilde{E}$ types, $\mathcal{K}$ is given by $X_{\delta}$, or a function thereof, defined as the image of any module with dimension $\delta$ under the Caldero--Chapoton map \cite{calderochapoton}, where $\delta$ gives the radical of the symmetrized Euler form \cite{crawleyboevey}. The values $\delta_i$ of $\delta$ at each vertex $i$ are given in Figure \ref{affinediagrams}.

We take the approach the authors of \cite{fordy,fordymarsh} used in the $\tilde{A}_{p,q}$ case. This is constructed by taking the $\tilde{A}$ diagram and orienting $p$ arrows clockwise and $q$ arrows anticlockwise. Due to the uniform nature of the mutations used, the frieze (\ref{generalfrieze}) can be written as the single recurrence
\begin{equation}\label{Arecurrence}
x_{n+N+1}x_n=x_{n+p}x_{n+q}+1,
\end{equation}
where $p+q=N+1$. The authors of \cite{fordy,fordymarsh} found the two periodic quantities
\[
J_n:=\frac{x_{n+2q}+x_n}{x_{n+q}}, \qquad \tilde{J}_n:=\frac{x_{n+2p}+x_n}{x_{n+p}}
\] 
which we display in the first row of Figure \ref{periodictable}. These immediately give the linear relations
\begin{equation}\label{Atypeperiodiclinearrelations}
x_{n+2q}-J_nx_{n+q}+x_n=0, \qquad x_{n+2p}-\tilde{J}_nx_{n+p}+x_n=0
\end{equation}
which are used to prove the $\tilde{A}$ case of (\ref{KSlinearrelation}) in \cite{fordy,fordymarsh}. In order to prove Theorem \ref{KSlinearrelation} for the $\tilde{D}$ and $\tilde{E}$ quivers we first find periodic quantities for their frieze sequences which give linear relations with periodic coefficients for the extending vertices. 
\begin{theorem}
The frieze sequences for the $\tilde{A}\tilde{D}\tilde{E}$ quivers have the periodic quantities given in Figure \ref{periodictable}. 
\end{theorem}
\begin{figure}
\begin{center}
\begin{tabular}{|c|c|c|}
\hline 
Quiver & Period & Quantity \\
\hline
& $p$ & $J_n$ \\
$\tilde{A}_{p,q}$ & & \\
& $q$ & $\tilde{J}_n$ \\
\hline
& $N-2$ & $J_n$ \\
$\tilde{D}_N$ & $2$ & $X^1_n/X^2_n$ \\
& $2$ & $X^N_n/X^{N+1}_n$ \\
\hline 
& $3$ & $J_n$ \\
$\tilde{E}_6$ & $3$ & $\tilde{J}_n$ \\
& $2$ & $K_n$ \\
\hline 
& $4$ & $J_n$ \\
$\tilde{E}_7$ & $3$ & $K_n$ \\
& $2$ & $\tilde{K}_n$ \\
\hline
 & $5$ & $J_n$ \\
$\tilde{E}_8$ & $3$ & $K_n$ \\
 & $2?$ & $\tilde{K}_n$\\
\hline
\end{tabular}
\end{center}\caption{Periodic quantities found for the $\tilde{A}\tilde{D}\tilde{E}$ quivers.}\label{periodictable}
\end{figure}
Here we find only the period $5$ and $3$ quantities, $J_n$ and $K_n$ respectively, for the $\tilde{E_8}$ quiver. The question mark is Conjecture \ref{period2conjecture} below. We remark the similarities between the widths of the tubes for the Auslander-Reiten quivers of $\tilde{A}\tilde{D}\tilde{E}$ type (see, for example, \cite{crawleyboevey}) and our Figure \ref{periodictable}. This is why we conjecture the missing period for $\tilde{E}_8$. 

For each quiver at least one of the quantities in Figure \ref{periodictable} gives a linear relation of the form (\ref{Atypeperiodiclinearrelations}) between the variables $X^k_n$ at only one vertex $k$. The rest of the quantities give linear relations that involve variables that live at different extending vertices. Focussing on the former, we derive relations that are remarkably similar to (\ref{Arecurrence}), two for $\tilde{E}_8$ (plus a third conjectured relation), and one for each of the remaining quivers. 
\begin{theorem}\label{similartoAtypetheorem}
For the $\tilde{D}_N$, with $N$ even, and $\tilde{E}$ quivers we have the following $\tilde{A}$ type recurrence:
\begin{equation}\label{recurrencesforextendingverticesequation}
X^k_{n+a+p}X^k_n=X^k_{n+a}X^k_{n+p}+\gamma_n^k.
\end{equation}
Here $k$ is an extending vertex and $\gamma^k_n$ is period $a$. For $\tilde{D}_N$ with $N$ odd we have instead
\begin{equation}\label{recurrencesforextendingverticesDodd}
X^1_{n+a+p}X^1_n=(\lambda_{n+1})^2X^1_{n+a}X^1_{n+p}+\gamma^1_n,
\end{equation}
and corresponding results for the other extending vertices, where $\lambda_{n+1}$ is period $2$ and $\gamma^1_n$ is period $a$. The values of $a$ and $p$ are given in Figure \ref{aroo}.
\begin{figure}
\begin{center}
\begin{tabular}{|c|c|c|}
\hline
Quiver & $a$ & $p$ \\
\hline\noalign{\smallskip}
$\tilde{D}_N \quad N \text{ even}$ & $1$ & $N-2$  \\
$\tilde{D}_N \quad N \text{ odd}$ & $1$ & $2N-4$ \\
$\tilde{E}_6$ & $3$ & $2$ \\
$\tilde{E}_7$ & $4$ & $3$\smallskip\smallskip \\
 & $6$ & $5$  \\
$\tilde{E}_8$ & $10$ & $3$  \\
 & $15?$ & $2?$  \\
\hline
\end{tabular}
\end{center}\caption{Values of $a$ and $p$ for the $\tilde{D}\tilde{E}$ quivers.}\label{aroo}
\end{figure}
\end{theorem}
Since a frieze is a system of recurrences, we can consider it a discrete dynamical system. In \cite{gsvomega} it was shown that a cluster algebra has a pre-symplectic form that is compatible with mutation. In our cases, however, this will be degenerate. We instead follow \cite{fordy} and project to a lower dimensional space with a reduced cluster map and a symplectic form.
We are interested in the integrability of this reduction, the definition of which we take from \cite{maeda,veselov}.
\begin{definition}\label{integrabilitydefinition}
Let $(M,\omega)$ be a $2m$-dimensional symplectic manifold and $\hat{\varphi}:M\rightarrow M$ a symplectic map. The map $\hat{\varphi}$ is called completely integrable if there exists $m$ functionally independent first integrals $\mathcal{I}_i$, $i=1,\ldots,m$, in involution, i.e $\{\mathcal{I}_i,\mathcal{I}_j\}=0$ for all $i,j$.
\end{definition}
\begin{theorem}\label{integrabilitytheorem}
The reduced cluster map for the $\tilde{D}_N$ quivers, where $N$ is odd, and for the $\tilde{E}$ type quivers is Liouville integrable. 
\end{theorem}
We provide an example of integrability for $\tilde{D}_6$ and leave the proof for general even $N$ open. \\\\
The structure of this paper is as follows: in Section $2$ we describe how the frieze can be obtained from cluster mutation and discuss our notion of integrability for cluster maps. In Section $3$ we explain our process for finding periodic quantities for the friezes and how we use these to construct linear relations. In Section $4$ we apply this to the $\tilde{D}$ and $\tilde{E}$ type quivers. In Section $5$ we prove integrability of the reduced cluster maps for these quivers.
\section{Background material on cluster algebras and integrability}\label{backgroundmaterial}
Here we briefly recall some notions from cluster algebras and discrete integrable systems. The following definitions largely follow \cite{fominintroduction} and \cite{fordy}. Here by a cluster algebra we mean a cluster algebra without coefficients or frozen variables. 

A quiver $Q$ is a directed graph where multiple edges are allowed. We disallow loops or 2-cycles. Quiver mutation $\mu_k$ at some vertex $k$ is defined in $3$ steps: 
\begin{enumerate}
\item For each length two path $i\rightarrow k \rightarrow j$ add a new arrow $i\rightarrow j$.
\item Reverse the direction of all arrows entering or exiting $k$.
\item Delete all 2-cycles that have appeared.
\end{enumerate}
The exchange matrix for a quiver $Q$ is the skew-symmetric matrix $B$ with entries $b_{ij}$ which is the number of arrows from $i$ to $j$, with $b_{ji}=-b_{ij}$. In addition to this, we will also attach a cluster variable $x_i$ at each vertex $i$. Cluster mutation at vertex $k$, also denoted $\mu_k$, fixes all variables $x_i$ with $i\neq k$ but transforms $x_k$ as follows:
\begin{equation}\label{clustermutationformula}
x'_k:=\mu_k(x_k):=\frac{1}{x_k}\left(\prod_{i\rightarrow k} x_j+\prod_{i\leftarrow k}x_i\right).
\end{equation}
Here the two products run over the vertices with arrows into and out of $k$ respectively. We consider $\mu_k$ as a mutation both of the quiver and of the cluster variables. One can see that if there is no edge between vertices $k$ and $k'$ then $\mu_k$ and $\mu_{k'}$ commute.

In the special case where the vertex $k$ is a sink or a source, the mutation simply reverses the arrows at $k$ and (\ref{clustermutationformula}) takes a simpler form, with one of the products being empty.
\subsection{Obtaining the frieze as a sequence of cluster mutations}\label{friezefromclusterseubsection}
To demonstrate the relationship between friezes and cluster algebras we recall how to obtain each frieze variable by cluster mutation, as described in \cite{kellerscherotzke}.

We first define initial cluster variables $X^k_0$ for $k=1,\ldots,N+1$ and take a labelling of the vertices such that $k$ is a sink in the full subquiver with vertices labelled greater than or equal to $k$, which is possible since the quiver is acyclic. The composition of mutations $\mu:=\mu_{N+1}\ldots\mu_2\mu_1$ then gives the frieze (\ref{generalfrieze}) for $n=1$, in the sense that $X^k_0$ becomes $X^k_{1}$ under this mutation. Moreover each vertex $k$ will be a sink when $\mu_k$ acts, so $\mu$ will flip each arrow twice, returning the original quiver. Further applications of $\mu$ will then give the frieze cluster variables $X^k_n$ for each $n$. We note that different labellings satisfying these conditions will give different friezes. 
\subsection{The bipartite belt as a frieze sequence}\label{bipartitebeltsubsection}
In the special case where the quiver is bipartite, one can take an orientation such that every vertex in either a sink or a source. We can then define the following compositions of mutations
\[
\mu_{\mathrm{sink}}:=\prod_{k \textrm{ a sink}}\mu_k, \qquad \mu_{\mathrm{source}}:=\prod_{k \textrm{ a source}}\mu_k,
\]
where the products are taken over vertices that are sinks (or sources) in the initial quiver. These products are well defined since no edge connects any pair of sinks (or sources). The composition $\mu:=\mu_{\mathrm{source}}\mu_{\mathrm{sink}}$ will, after relabelling of vertices, satisfy the conditions of the previous subsection, and in this case the frieze relations (\ref{generalfrieze}) become simpler: for each $k=1,\ldots,N+1$ we have
\begin{equation}\label{friezemutations}
\begin{aligned}
X_{n+1}^kX^k_n=\left(1+\prod_{j\rightarrow k}X^j_n\right) & \qquad \text{ if } k \text{ is a sink, } \\
X_{n+1}^kX^k_n=\left(1+\prod_{j\leftarrow k}X^j_{n+1}\right) & \qquad \text{ if } k \text{ is a source. } 
\end{aligned}
\end{equation}
The cluster variables constructed in this way form the bipartite belt of the cluster algebra associated to the quiver, defined in \cite{clustersiv}. 
\subsection{Integrability of dynamical systems from periodic quivers}\label{dynamicalsystemsubsection}
The pre-symplectic form associated with the exchange matrix $B$ is
\[
\omega:=\sum_{i<j}\frac{b_{ij}}{x_ix_j}dx_i\wedge dx_j
\]
and was shown in \cite{gsvomega} to be compatible with cluster mutation, in the sense that it transforms as
\[
\mu_k^\ast(\omega)=\sum_{i<j}\frac{b'_{ij}}{x'_ix'_j}dx'_i\wedge dx'_j,
\] 
where the $x_i$ are generic cluster variables, $b_{ij}$ are the entries of the exchange matrix and the prime denotes the images of these under the mutation $\mu_k$. In our work the cluster variables will be only those that are part of the frieze sequence, so the $x_i$ will all be of the form $X^i_n$.  
\begin{lemma}
The cluster map $\varphi$ associated with a frieze sequence preserves the symplectic form, i.e. $\varphi^*\omega=\omega$.
\begin{proof}
As in the $\tilde{A}$ case discussed in \cite{fordyhone} the cluster map acts as
\[
\varphi^*\omega=
\sum_{1\leq i<j\leq N+1}\frac{b_{ij}}{X^i_{n+1}X^j_{n+1}}dX^i_{n+1}\wedge dX^j_{n+1}.
\]
This, however, may be decomposed into $\varphi=\varphi_{\text{source}}\circ\varphi_{\text{sink}}$ where $\varphi_{\text{sink}}$ and $\varphi_{\text{source}}$ are the maps given by $\mu_{\mathrm{sink}}$ and $\mu_{\mathrm{source}}$ respectively. We have
\[
\varphi_{\text{sink}}^*(\omega)=
\sum_{\substack{ i \text { a sink }\\
j \text{ a source }}}
\frac{b_{ij}}{X^i_{n+1}X^j_{n}}dX^i_{n+1}\wedge dX^j_{n}.
\]
Taking $i$ to be a sink the mutation relation is 
\[
X^i_{n+1}X^i_n=1+\prod_{j\rightarrow  i}X^j_n
\]
giving
\[
dX^i_{n+1}=-\frac{X^i_{n+1}}{X^i_n}dX^i_n+\frac{X^i_{n+1}X^i_n-1}{X^i_n}\sum_{j\rightarrow i}\frac{dX^j_n}{X^j_n}.
\]
We examine the following sum:
\[
\sum_{j\rightarrow i}\frac{dX^i_{n+1}\wedge dX^j_n}{X^i_{n+1}X^j_n}
=
-\sum_{j\rightarrow i}\frac{dX^i_{n}\wedge dX^j_n}{X^i_{n}X^j_n}
+
\sum_{j\rightarrow i}\frac{X^i_{n+1}X^i_n-1}{X^i_n}\sum_{k\rightarrow i}\frac{dX^k_{n}\wedge dX^j_n}{X^k_{n}X^j_n}
\]
\[
=-\sum_{j\rightarrow i}\frac{dX^i_{n}\wedge dX^j_n}{X^i_{n}X^j_n}
+
\frac{X^i_{n+1}X^i_n-1}{X^i_n}\sum_{k,j\--i}\frac{dX^k_{n}\wedge dX^j_n}{X^k_{n}X^j_n}
\]
and see that the final sum is over every $k,j$ connected to $i$, hence each $\frac{dX^k_{n}\wedge dX^j_n}{X^k_{n}X^j_n}$ is counted twice, with opposite sign, giving zero. Hence
\[
\sum_{j\rightarrow i}\frac{dX^i_{n+1}\wedge dX^j_n}{X^i_{n+1}X^j_n}
=
-\sum_{j\rightarrow i}\frac{dX^i_{n}\wedge dX^j_n}{X^i_{n}X^j_n}
\]
so $\varphi_{\text{sink}}^*(\omega)=-\omega$. A similar proof shows that $\varphi_{\text{source}}^*(\omega)=-\omega$ so the decomposition $\varphi=\varphi_{\text{source}}\circ\varphi_{\text{sink}}$ gives the result.
\end{proof}
\end{lemma}
In the cases we shall encounter the determinant of $B$ will be zero, so $\omega$ will be degenerate. To remedy this we project to a lower dimensional space where $\varphi$ reduces to a symplectic map. For each $n$ we take reduced coordinates $\{y^j_n\}_{j=1,\ldots,2m}$, where $\rank(B)=2m$ and
\[
y^j_n:=\prod_i (X^i_n)^{v_{j,i}}.
\]
Here $v_{j,i}$ is the $i$th component of $v_j$ and the vectors $v_1,\ldots,v_{2m}$ are a $\mathbb{Z}$-basis for $\image{B}$. We define this projections by
\[
\pi:\{X^1_n,X^2_n,\ldots ,X^{N+1}_n\}\rightarrow\{y^1_n,y^2_n,\ldots,y^{2m}_n\}.
\] 
As in Lemma 2.7 of \cite{fordy} we have that for each $n$ the set $\{y^j_n\}_{j=1,\ldots,2m}$ is a full set of invariants for the group of scaling transformations defined by 
\begin{equation}\label{lambdascaling}
(X^1_n,X^2_n\ldots)\mapsto (\lambda^{u_1}X^1_n,\lambda^{u_{2}}X^2_n,\ldots,)
\end{equation}
for any $\lambda\in \mathbb{C}^{\star}$ and $u=(u_1,\ldots,u_N)\in \kernel B$, and that for each $n$ the two form 
\begin{equation}\label{symplectictwoform}
\hat{\omega}=\sum_{i<j}\frac{\hat{b}_{ij}}{y^i_ny^j_n}dy^i_n \wedge dy^j_n,
\end{equation}
is symplectic with $\pi^*\circ\hat{\omega}=\omega$. Here the coefficients $\hat{b}_{ij}$ are given by the entries of $\hat{B}$ which satisfies 
\[
A^{-T}BA^{-1}=
\begin{pmatrix}
\hat{B} & 0 \\
0 & 0 
\end{pmatrix}
\]
where the matrix $A$ has the vectors $v_1,\ldots, v_{N+1}$ as rows, and the sets of vectors $\{v_1,\ldots,v_{2m}\}$ and $\{v_{2m+1},\ldots, v_{N+1}\}$ are bases for $\image{B}$ and $\kernel{B}$ respectively. The form of $\hat{B}$ and the reduced coordinates depend on the choice of basis for the image and kernel of $B$. The following theorem is the analogue of Theorem 2.6 in \cite{fordy}.
\begin{theorem}
For the $\tilde{D}\tilde{E}$ type quivers the cluster map $\varphi$ reduces to a map $\hat{\varphi}$:
\begin{equation}\label{reducedclustermap}
\hat\varphi:
\begin{pmatrix}
y^1_n \\
y^2_n \\
\vdots \\
y^{2m}_n
\end{pmatrix}
\mapsto
\begin{pmatrix}
y^1_{n+1} \\
y^2_{n+1} \\
\vdots \\
y^{2m}_{n+1}
\end{pmatrix}
\end{equation} 
such that $\pi\circ\varphi=\hat{\varphi}\circ\pi$. The map $\hat{\varphi}$ is birational and symplectic with respect to the two form (\ref{symplectictwoform}). 
\begin{proof}
The property $\pi\circ\varphi=\hat{\varphi}\circ\pi$ and birationality of the maps $\hat{\varphi}$ is shown as they are constructed, case by case, in Section \ref{integrabilitysection}. That $\hat{\varphi}$ is symplectic is the same as in \cite{fordy} Lemma 2.9.
\end{proof}
\end{theorem}
Since $\hat{\omega}$ is symplectic, it defines a nondegenerate Poisson bracket:
\begin{equation}\label{logcanonicalbracket}
\{y^i_n,y^j_n\}=c_{ij}y^i_ny^j_n, \quad \text { where } \quad C=(c_{ij}):=\hat{B}^{-1}.
\end{equation}
Our aim is to prove Liouville integrability of the map (\ref{reducedclustermap}), as in Definition \ref{integrabilitydefinition}, by finding enough involutive first integrals with respect to this bracket.
\section{Method of proof for the $\tilde{D}$ and $\tilde{E}$ quivers}\label{methodofproofsection}
Here we'll discuss our general process for finding periodic quantities for the frieze sequences. We show how to find constant coefficient linear relations from these. 

We orient the $\tilde{D}$ and $\tilde{E}$ quivers such that every vertex is either a sink or a source. This allows us to construct the bipartite belt frieze described in Subsection \ref{bipartitebeltsubsection}.

The key to our results, for $\tilde{D}$ type, is to write each of the mutation relations (\ref{friezemutations}) as a $2\times 2$ matrix with determinant $-1$:
\begin{equation}\label{clustermutationdeterminantform}
-1=
\begin{vmatrix}
\star & X^k_{n+1}  \\
X^k_n &\star
\end{vmatrix},
\end{equation}
where the product of the stars will give the product inside the brackets in (\ref{friezemutations}). For the $\tilde{E}$ quivers we have taken the determinant in (\ref{clustermutationdeterminantform}) to be $+1$, but the procedure is equivalent, up to swapping matrix columns. We construct a $3\times 3$ matrix such that each connected $2\times 2$ determinant inside is of the form (\ref{clustermutationdeterminantform}). By the following lemma, part of the Dodgson condensation algorithm \cite{dodgson}, the $3\times 3$ matrix will have determinant zero. 
\begin{lemma}\label{dodgsonlemma}
For general $x_{ij}$ we have 
\begin{equation}\label{dodgson}
\begin{vmatrix}
x_{0,0} & x_{0,1} & x_{0,2} \\
x_{1,0} & x_{1,1} & x_{1,2} \\
x_{2,0} & x_{2,1} & x_{2,2}  
\end{vmatrix}
x_{1,1}
=\begin{vmatrix}
x_{0,0} & x_{0,1} \\
x_{1,0} & x_{1,1}
\end{vmatrix}
\begin{vmatrix}
x_{1,1} & x_{1,2} \\
x_{2,1} & x_{2,2}
\end{vmatrix}
-
\begin{vmatrix}
x_{0,1} & x_{0,2} \\
x_{1,1} & x_{1,2} 
\end{vmatrix}
\begin{vmatrix}
x_{1,0} & x_{1,1} \\
x_{2,0} & x_{2,1}
\end{vmatrix}.
\end{equation}
\end{lemma}
In our cases the $3\times 3$ matrix will have a kernel vector of the form $(1,-\alpha,1)^T$. We then add rows to the bottom (or top) of this matrix that annihilate this vector. Our main method for doing this will be such that, using (\ref{clustermutationdeterminantform}), the new row will create two new connected $2 \times 2$ matrices with equal determinant. By the following lemma this newly formed $4\times 3$ matrix will have the same kernel vector.
\begin{lemma}\label{kernelpreservation}
Let $M$ be a $3 \times 3$ matrix 
\[
M:=
\begin{pmatrix}
x_{0,0} & x_{0,1} & x_{0,2} \\
x_{1,0} & x_{1,1} & x_{1,2} \\
x_{2,0} & x_{2,1} & x_{2,2}  
\end{pmatrix}
\]
such that $(1,-\alpha,1)^T$ is in the kernel of the first two rows. We define the determinants $\delta_{i,j}$ by
\[
\delta_{i,j}:=
\begin{vmatrix}
x_{i,j} & x_{i,j+1} \\
x_{i+1,j} & x_{i+1,j+1}
\end{vmatrix}.
\] 
If $x_{1,1}\neq 0$ and $\delta_{1,0}=\delta_{1,1}\neq 0$ then $(1,-\alpha,1)^T$ is in the kernel for all three rows of $M$.
\begin{proof}
By (\ref{dodgson}), and since $\delta_{1,0}=\delta_{1,1}\neq0$, we have
\[
|M|=0 \Leftrightarrow \delta_{0,0}\delta_{1,1}-\delta_{1,0}\delta_{0,1}=0 \Leftrightarrow \delta_{0,0}-\delta_{0,1}=0.
\]
Now the matrix 
\[
\begin{pmatrix}
x_{0,0} & x_{0,1} & x_{0,2} \\
x_{1,0} & x_{1,1} & x_{1,2} \\
-1 & 0 & 1 
\end{pmatrix}
\]
has a kernel vector $(1,-\alpha,1)^T$ and determinant zero. We expand this determinant along the bottom row to see  that $\delta_{0,0}=\delta_{0,1}$. Hence $|M|=$ is zero. Taking a kernel vector $(A,B,1)^T$ for $M$ we have $A=\frac{\delta_{1,1}}{\delta_{1,0}}=1$ by Cramer's rule and $B=-\alpha$.
\end{proof}
\end{lemma}
After adding one row our $4\times 3$ matrix will be of the form
\[
\begin{pmatrix}
x_{0,0} & x_{0,1} & x_{0,2} \\
x_{1,0} & x_{1,1} & x_{1,2} \\
x_{2,0} & x_{2,1} & x_{2,2} \\
x_{3,0} & x_{3,1} & x_{3,2}
\end{pmatrix}
\] 
so that the bottom $3 \times 3$ matrix satisfies the conditions of Lemma \ref{kernelpreservation}, hence the $4\times 3$ matrix has $(1,-\alpha,1)^T$ as a kernel vector. As mentioned in \ref{clustermutationdeterminantform} for the $x_{i,j}$ we will take cluster variables $X^k_n$, or functions thereof. We continue adding rows in this way (except for a few exceptional rows in the $\tilde{D}$ cases) until we find that either $\alpha$ is periodic, in the $\tilde{D}$ cases, or one of the matrix entries is periodic, in the $\tilde{E}$ cases, with respect to the shift $\varphi^*$.
\subsection{Linear relations with fixed coefficients}
The periodic quantities we find immediately give linear relations with periodic coefficients, which in turn give linear relations with constant coefficients. In each case we'll have $\Psi_nL_n=\Psi_{n+q}$ for some integer $q$, where $\Psi_n$ is a $2 \times 2$ matrix with cluster variable entries and $L_n$ is a $2 \times 2$ matrix with periodic elements, i.e. $L_{n+p}=L_n$ for some $p$. Taking $m:=\lcm(p,q)$ we have
\begin{equation}\label{generalM}
\Psi_nM_n=\Psi_{n+m}, \qquad M_n:=L_nL_{n+q}\ldots L_{n+m-q}.
\end{equation}
\begin{theorem}\label{CayleyHamilton}
The matrix $\Psi_n$ satisfies
\begin{equation}\label{generalpsi}
\Psi_{n+2m}-\trace{(M_n)}\Psi_{n+m}+\det(M_n)\Psi_n=0,
\end{equation}
where $\trace{(M_n)}$ is period $q$.
\begin{proof}
Applying the Cayley-Hamilton theorem to $M_n$ gives 
\[
M_n^2-\trace(M_n)M_n+\det(M_n)I=0.
\] 
Multiplying by $\Psi_n$ from the left gives the result.
\end{proof}
\end{theorem}
In each of the cases we deal with, however, we find that $\mathcal{K}:=\trace(M_n)$ is invariant (under shifts in $n$) and that $\det(M_n)=1$. Hence (\ref{generalpsi}) reduces to $\Psi_{n+2m}-\mathcal{K}\Psi_{n+m}+\Psi_n=0$ and any entry of this matrix equation will give us constant coefficient linear relations between cluster variables.
\section{Linear relations in the $\tilde{D}$ and $\tilde{E}$ Quivers}
Here we apply the ideas of the previous section to the $\tilde{D}$ and $\tilde{E}$ quivers. For $\tilde{D}_N$ we deal with different parities separately. As well as those stated in the introduction we have sporadic appearances of other linear relations.
\subsection{Quivers of $\tilde{D}$ type.}
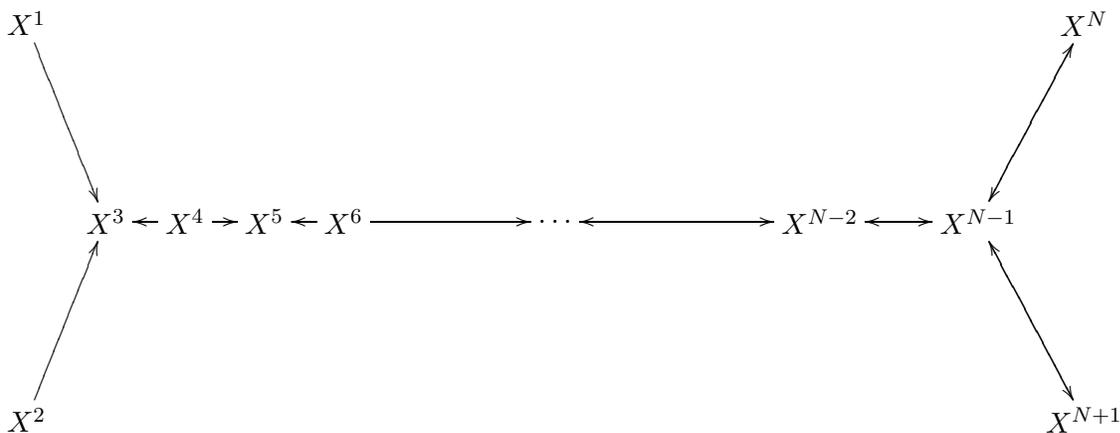
\begin{figure}[h]
\begin{xy} 0;<1pt,0pt>:<0pt,-1pt>:: 
(0,0) *+{X^1} ="X^1",
(0,150) *+{X^2} ="X^2",
(60,75) *+{X^4} ="X^4",
(30,75) *+{X^3} ="X^3",
(90,75) *+{X^5} ="X^5",
(120,75) *+{X^6} ="X^6",
(200,75) *+{\ldots} ="\ldots",
(300,75) *+{X^{N-2}} ="X^{N-2}",
(360,75) *+{X^{N-1}} ="X^{N-1}",
(400,0) *+{X^N} ="X^N",
(400,150) *+{X^{N+1}} ="X^{N+1}",
"X^1", {\ar"X^3"},
"X^2", {\ar"X^3"},
"X^4", {\ar"X^3"},
"X^4", {\ar"X^5"},
"X^6", {\ar"X^5"},
"X^6", {\ar"\ldots"},
"\ldots", {\ar"X^{N-2}"},
"X^{N-2}", {\ar"\ldots"},
"X^{N-2}", {\ar"X^{N-1}"},
"X^{N-1}", {\ar"X^{N-2}"},
"X^N", {\ar"X^{N-1}"},
"X^{N-1}", {\ar"X^N"},
"X^{N+1}", {\ar"X^{N-1}"},
"X^{N-1}", {\ar"X^{N+1}"},
\end{xy}
\caption{The $\tilde{D}_N$ Quiver.}\label{Dquiver}
\end{figure}
We'll orient the $\tilde{D}_N$ diagram as in Figure \ref{Dquiver}. The two sided arrows are used on the right-hand side to indicate that the direction depends on the parity of $N$, i.e. for a given $N$ one of the arrowheads needs to be deleted on each double arrow, ensuring that each vertex is a sink or source. We deal with the case $N=4$ separately, in Section \ref{N=4case}. 

We use the sequence of mutations described in Section \ref{methodofproofsection}. For example, mutating at vertex $3$ gives
\begin{equation}\label{clusterrelation}
X^3_{n+1}=\frac{1}{X^3_n}(1+X^1_nX^2_nX^4_n)
\end{equation}
and then at vertex $1$ 
\[
X^1_{n+1}=\frac{1}{X^1_n}\left(1+X^3_{n+1}\right)=\frac{1}{X^1_n}\left(1+\frac{1}{X^3_n}(1+X^1_nX^2_nX^4_n)\right).
\]
In general the relations are
\begin{align}
X^1_{n+1}&=\frac{1}{X^1_n}(1+X^3_{n+1}), & X^2_{n+1}=&\frac{1}{X^2_n}(1+X^3_{n+1}), &X^3_{n+1}=&\frac{1}{X^3_{n}}(1+X^1_nX^2_nX^4_n), \nonumber\\
X^i_{n+1}&=\frac{1}{X^i_{n}}(1+X^{i-1}_{n+1}X^{i+1}_{n+1}),& &
3<i<N-1 &  &\text{ for } i \text{ even,} \nonumber \\
X^i_{n+1}&=\frac{1}{X^i_{n}}(1+X^{i-1}_{n}X^{i+1}_{n}),& &
3<i<N-1 &  &\text{ for } i \text{ odd.} \label{Drelations1}
\end{align}
Towards the right end of the diagram we need to distinguish two cases, depending on the parity of $N$. When $N$ is even we have
\begin{align}\label{Drelations2}
X^{N-1}_{n+1}&=\frac{1}{X^{N-1}_n}(1+X^{N-2}_nX^N_nX^{N+1}_n), \nonumber\\
X^N_{n+1}&=\frac{1}{X^N_{n}}(1+X^{N-1}_{n+1}),&
X^{N+1}_{n+1}=\frac{1}{X^{N+1}_{n}}(1+X^{N-1}_{n+1}).
\end{align}
and when $N$ is odd we have
\begin{align}\label{Drelations3}
X^{N-1}_{n+1}&=\frac{1}{X^{N-1}_n}(1+X^{N-2}_{n+1}X^N_{n+1}X^{N+1}_{n+1}), \hspace{10 mm} \nonumber \\
X^N_{n+1}&=\frac{1}{X^N_{n}}(1+X^{N-1}_{n}),&
X^{N+1}_{n+1}=\frac{1}{X^{N+1}_{n}}(1+X^{N-1}_{n}).
\end{align}
\begin{lemma}\label{Dtypeperiod2}
We note that from (\ref{Drelations1}) at the extending vertices $1$ and $2$ we have
\[
\frac{X^1_{n+1}}{X^2_{n+1}}=\frac{X^2_{n}}{X^1_{n}}=\ldots =\left(\frac{X^1_{0}}{X^2_{0}}\right)^{(-1)^n},
\]  
as well as a corresponding result for vertices $N$ and $N+1$ from Equations \ref{Drelations2} and \ref{Drelations3}. In particular $X^1_{n}/X^2_n$ and $X^{N+1}_{n}/X^N_n$ are period $2$.
\end{lemma}
The $3\times 3$ matrix from Lemma \ref{kernelpreservation} that we'll use is the following:
\begin{lemma} 
The matrix
\[
M:=\begin{pmatrix}
X^1_{n+1}X^2_{n+1} & X^3_{n+1} & X^4_{n}  \\
X^3_{n+2} & X^4_{n+1} & X^5_{n+1}  \\
X^4_{n+2} & X^5_{n+2} & X^6_{n+1}  \\
\end{pmatrix}
\]
has kernel vector $(1,-J_n,1)^T$, where
\[
J_n:=\frac{X^1_{n+1}X^2_{n+1}+X^4_n}{X^3_{n+1}}.
\]
\begin{proof}
Each of the cluster relations, for example (\ref{clusterrelation}), can be written as
\[
-1=
\begin{vmatrix}
\star & \star \\
\star & \star
\end{vmatrix}
\]
for appropriate $2\times 2$ matrices, so each connected $2 \times 2$ determinant inside $M$ is equal to $-1$. We see from (\ref{Drelations1}) that
\[
\left(1,-\frac{X^1_{n+1}X^2_{n+1}+X^4_n}{X^3_{n+1}},1\right)^T
\]
is in the kernel of the first two rows of $M$. Hence by Lemma \ref{kernelpreservation} this vector is in the kernel of $M$. 
\end{proof}
\end{lemma}
We can extend $M$ to involve the variables at each vertex, ensuring that each connected $2\times 2$ matrix has determinant $-1$, as described in Section \ref{methodofproofsection}. By repeated applications of Lemma \ref{kernelpreservation} the vector $(1,-J_n,1)$ will be in the kernel of each newly added row. For example, after the first application we have
\begin{equation*}
M:=\begin{pmatrix}
X^1_{n+1}X^2_{n+1} & X^3_{n+1} & X^4_{n}  \\
X^3_{n+2} & X^4_{n+1} & X^5_{n+1}  \\
X^4_{n+2} & X^5_{n+2} & X^6_{n+1}  \\
X^5_{n+3} & X^6_{n+2} & X^7_{n+2} 
\end{pmatrix},
\end{equation*}
with the newly formed connected $2 \times 2$ matrices
\[
\begin{pmatrix}
X^4_{n+2} & X^5_{n+2}  \\
X^5_{n+3} & X^6_{n+2}
\end{pmatrix},
\qquad 
\begin{pmatrix}
X^5_{n+2} & X^6_{n+1}  \\
X^6_{n+2} & X^7_{n+2} 
\end{pmatrix}
\]
each having determinant $-1$, which follows from the mutation relations (\ref{Drelations1}). The kernel vector $(1,-J_n,1)$ for this $4 \times 3$ matrix gives us a new equation
\[
\begin{pmatrix}
X^5_{n+3} & X^6_{n+2} & X^7_{n+2}
\end{pmatrix}
\begin{pmatrix}
1 \\
-J_n \\
1
\end{pmatrix}
=0.
\]
This particular relation is not so important, we simply state it to stress that we arrive at a new equation for each row we add. We can use the Equations \ref{Drelations1} to add $N-6$ rows to $M$ in the same way:
\begin{equation*}\label{matrix}
M:=\begin{pmatrix}
X^1_{n+1}X^2_{n+1} & X^3_{n+1} & X^4_{n}  \\
X^3_{n+2} & X^4_{n+1} & X^5_{n+1}  \\
X^4_{n+2} & X^5_{n+2} & X^6_{n+1}  \\
X^5_{n+3} & X^6_{n+2} & X^7_{n+2}  \\
\vdots & \vdots & \vdots  \\
X^{N-3}_{n+l} & X^{N-2}_{n+l-h} & X^{N-1}_{n+l-1} \\
X^{N-2}_{n+l+1-h} & X^{N-1}_{n+l} & X^{N}_{n+l-h}X^{N+1}_{n+l-h} \\
\end{pmatrix}.
\end{equation*} 
Here the subindex in the leftmost column increases by $1$ only at each sink, so $l$ is $1$ greater than the number of sinks from $X^3$ to $X^{N-3}$ inclusive. The value of $h$ is $0$ if $N$ is odd and $1$ if $N$ is even. To proceed we can use the identities:
\begin{equation}\label{Dtypeidentities}
J_n=\frac{X^1_{n+1}+X^1_{n-1}}{X^2_n}=\frac{X^1_{n+1}X^2_{n+1}+X^4_n}{X^3_{n+1}}=\frac{X^1_{n-1}X^2_{n-1}+X^4_n}{X^3_n},
\end{equation} 
which follow from some simple manipulation of (\ref{Drelations1}), to add two rows (shown in black) to the top of the following matrix while preserving our kernel vector. Note that these new connected $2\times 2$ minors don't have determinant $-1$. This gives us
\[
M:=
\begin{pmatrix}
X_n^4 & X_n^3 & X_{n-1}^1X_{n-1}^2 \\
X^1_{n+1} & X_n^2 & X^1_{n-1} \\
\color{grey}X^1_{n+1}X^2_{n+1} & \color{grey}X^3_{n+1} & \color{grey}X^4_{n}  \\
\vdots & \vdots & \vdots  \\
\color{grey}X^{N-2}_{n+l+1-h} & \color{grey}X^{N-1}_{n+l} & \color{grey}X^{N}_{n+l-h}X^{N+1}_{n+l-h} \\
\end{pmatrix}.
\] 
In order to extend the bottom in this way we'll need to use similar identities to (\ref{Dtypeidentities}) for the vertices at the right end of the diagram, but these depend on the parity of $N$, so now we split the problem in to two cases. 

Before we do this we note that for small $N$ the above process may be not defined. We deal with the case $N=4$ in Section \ref{N=4case}. For $N=5$ one should start with the matrix 
\[
\begin{pmatrix}
X^1_{n+1}X^2_{n+1} & X^3_{n+1} & X^4_n \\
X^3_{n+2} & X^4_{n+1} & X^5_{n+1}X^6_{n+1}
\end{pmatrix}
\]
and use (\ref{Dtypeidentities}) to add rows above. Rows below should then be added as described in Section \ref{oddNsubsection}. For $N=6$ add rows above the matrix
\[
\begin{pmatrix}
X^1_{n+1}X^2_{n+1} & X^3_{n+1} & X^4_n \\
X^3_{n+2} & X^4_{n+1} & X^5_{n+1} \\
X^4_{n+2} & X^5_{n+2} & X^6_{n+1}X^7_{n+1}
\end{pmatrix}
\]
with (\ref{Dtypeidentities}) and add rows below as described in Section \ref{evenNsubsection}.
\subsubsection{The even $N$ case with $N>4$.}\label{evenNsubsection}
Here the vertices $N$ and $N+1$ are sources, $l=\frac{N}{2}-1$ and $h=1$. We now establish why $J_n$ is of such importance.
\begin{lemma}\label{periodiclemma} The quantity of (\ref{Dtypeidentities}), $J_n=\frac{X^1_{n+1}+X^1_{n-1}}{X^2_n}$, is periodic with period $2l=N-2$
\begin{proof}
The diagram is symmetric about the centre vertex so the relations
 \begin{equation}\label{evenDtyperelations}
\frac{X^{N+1}_{n+1}+X^{N+1}_{n-1}}{X^N_n}=\frac{X^{N+1}_{n+1}X^N_{n+1}+X^{N-2}_n}{X^{N-1}_{n+1}}=\frac{X^{N+1}_{n-1}X^N_{n-1}+X^{N-2}_n}{X^{N-1}_n}
\end{equation}
mirror (\ref{Dtypeidentities}). We use these to extend $M$:
\[
\begin{pmatrix}
\color{grey}X_n^4 & \color{grey}X_n^3 & \color{grey}X_{n-1}^1X_{n-1}^2 \\
\color{grey}X^1_{n+1} & \color{grey}X_n^2 & \color{grey}X^1_{n-1} \\
\vdots & \vdots & \vdots  \\
\color{grey}X^{N-2}_{n+l} & \color{grey}X^{N-1}_{n+l} & \color{grey}X^{N}_{n+l-1}X^{N+1}_{n+l-1} \\
X_{n+l+1}^{N+1} & X_{n+l}^N & X_{n+l-1}^{N+1} \\
X_{n+l+1}^{N+1}X^N_{n+l+1} & X_{n+l+1}^{N-1} & X_{n+l}^{N-2}
\end{pmatrix}.
\]
Now the bottom row here is the counterpart to the top row of (\ref{matrix}) from the opposite end of the quiver. Hence we can add extra rows below in exactly the same way in which we started:
\[
\begin{pmatrix}
\color{grey}X_n^4 & \color{grey}X_n^3 & \color{grey}X_{n-1}^1X_{n-1}^2 \\
\color{grey}X^1_{n+1} & \color{grey}X_n^2 & \color{grey}X^1_{n-1} \\
\color{grey}X^1_{n+1}X^2_{n+1} & \color{grey}X^3_{n+1} & \color{grey}X^4_{n}  \\
\color{grey} X^3_{n+2} & \color{grey}X^4_{n+1} & \color{grey}X^5_{n+1}  \\
\vdots & \vdots & \vdots  \\
\color{grey}X_{n+l+1}^{N+1} & \color{grey}X_{n+l}^N & \color{grey}X_{n+l-1}^{N+1} \\
\color{grey}X_{n+l+1}^{N+1}X^N_{n+l+1} & \color{grey}X_{n+l+1}^{N-1} & \color{grey}X_{n+l}^{N-2} \\
X^{N-1}_{n+l+2} & X^{N-2}_{n+l+1} & X^{N-3}_{n+l+1} \\
X^{N-2}_{n+l+2} & X^{N-3}_{n+l+2} & X^{N-4}_{n+l+1} \\
\vdots & \vdots & \vdots \\
X^5_{n+2l} & X^4_{n+2l-1} & X^3_{n+2l-1} \\
X^4_{n+2l} & X^3_{n+2l} & X^2_{n+2l-1}X^1_{n+2l-1} \\
X^1_{n+2l+1} & X^2_{n+2l} & X^1_{n+2l-1}
\end{pmatrix},
\] 
where the final row follows from (\ref{Dtypeidentities}). Comparing the second and last rows we can see that 
\[
J_n=\frac{X^1_{n+1}+X^1_{n-1}}{X^2_n}=\frac{X^1_{n+2l+1}+X^1_{n+2l-1}}{X^2_{n+2l}}=J_{n+2l}.
\]
\end{proof}
\end{lemma}
\begin{theorem}\label{evenDlinearrelationtheorem}
For even $N>4$ the constant coefficient linear relation for the variables at the extending vertices is
\begin{equation}\label{evenDlinearrelation}
X^k_{n+2N-4}-\mathcal{K}X^k_{n+N-2}+X^k_n=0,
\end{equation} 
by which we mean that $k\in\{1,2,N,N+1\}$ and $\mathcal{K}$ is invariant.
\begin{proof}
Defining 
\[
\Psi_n:=
\begin{pmatrix}
X^1_n & X^2_{n+1} \\
X^1_{n+N-2} & X^2_{n+N-1}
\end{pmatrix},
\quad
\tilde{L}_n:=
\begin{pmatrix}
0 & -1 \\
1 & J_n
\end{pmatrix},
\quad
\text{ then }
\quad
\Psi_n\tilde{L}_n=
\begin{pmatrix}
X^2_{n+1} & X^1_{n+2} \\
X^2_{n+N-1} & X^1_{n+N}
\end{pmatrix},
\]
so $\Psi_n\tilde{L}_n\tilde{L}_{n+1}=\Psi_{n+2}$. Calling $L_n:=\tilde{L}_n\tilde{L}_{n+1}$ allows us to apply Theorem \ref{CayleyHamilton} with $q=2$ and $p=m=N-2$. Explicitly $M_n=L_nL_{n+2}\ldots L_{n+N-4}$. Since $\tilde{L}_{n+N-2}=\tilde{L}_n$ and the trace is fixed under cyclic permutations we have that $\trace(M_n)=\mathcal{K}$ is invariant. We also have $|M_n|=1$. The theorem gives $\Psi_{n+2N-4}-\mathcal{K}\Psi_{n+N-2}+\Psi_n=0$ and the top left entry of this matrix equation gives the linear relation (\ref{evenDlinearrelation}) for $k=1$. By symmetry this holds for the other extending vertices.
\end{proof}
\end{theorem} 
We now give the following result, one of the cases in Theorem \ref{similartoAtypetheorem}.
\begin{corollary}\label{evenNcorollary}
For even $N>4$ the $X^1_n$ variables satisfy the $\tilde{A}$ type recurrence
\[
X^1_{n+N-1}X^1_n=X^1_{n+N-2}X^1_{n+1}+\gamma^1,
\]
where $\gamma^1$ is invariant. 
\begin{proof}
Taking determinants in $\Psi_n\tilde{L}_n=\Psi_{n+1}$ we have that
\[
X^1_nX^2_{n+N-1}-X^2_{n+1}X^1_{n+N-2}=X^2_{n+1}X^1_{n+N}-X^1_{n+2}X^2_{n+N-1},
\]
and by Lemma \ref{Dtypeperiod2} we can write $X^2_n$ as $X^1_n\lambda_n$ where $\lambda_n:=\left(\frac{X^2_0}{X^1_0}\right)^{(-1)^n}$ which is period $2$. When replacing each $X^2_n$, each $\lambda$ will appear with the same subscript and we cancel them to give the result.
\end{proof}
\end{corollary}
\subsubsection{The case $N=4$.}\label{N=4case}
In this case the mutation relations are given by
\begin{align}
X^i_{n+1}&=\frac{1}{X^i_n}(1+X^3_{n+1}), & \textrm{for } i=1,2,4,5, \nonumber\\
X^3_{n+1}&=\frac{1}{X^3_n}(1+X^1_nX^2_nX^4_nX^5_n). \nonumber
\end{align} 
The analogue of (\ref{Dtypeidentities}) is
\[
J_n:=\frac{X^1_{n+1}+X^1_{n-1}}{X^2_n}=\frac{X^1_{n+1}X^2_{n+1}+X^4_nX^5_n}{X^3_{n+1}}=\frac{X^1_{n-1}X^2_{n-1}+X^4_nX^5_n}{X^3_n}
\]
and the analogous result
\[
\sigma(J_n)=\frac{X^4_{n+1}+X^4_{n-1}}{X^5_n}=\frac{X^4_{n+1}X^5_{n+1}+X^1_nX^2_n}{X^3_{n+1}}=\frac{X^4_{n-1}X^5_{n-1}+X^1_nX^2_n}{X^3_n}
\]
is given by the permutation $\sigma:=(14)(25)$. A simple calculation shows that $J_n=\sigma(J_{n-1})$ hence $J_n$ is period $N-2=2$. Since the expression $J_n=\frac{X^1_{n+1}+X^1_{n-1}}{X^2_n}$ is the same as the other cases, Theorem \ref{evenDlinearrelationtheorem} and Corollary \ref{evenNcorollary} hold for $N=4$.
\subsubsection{The odd $N$ case.}\label{oddNsubsection}
Here vertices $N$ and $N+1$ are sinks, $l=\frac{N-3}{2}$, and $h=0$. We can start with the matrix
\[
\begin{pmatrix}
    X_n^4 & X_n^3 & X_{n-1}^1X_{n-1}^2 \\
    X^1_{n+1} & X_n^2 & X^1_{n-1} \\
    X^1_{n+1}X^2_{n+1} & X^3_{n+1} & X^4_{n}  \\
    X^3_{n+2} & X^4_{n+1} & X^5_{n+1}  \\
   X^4_{n+2} & X^5_{n+2} & X^6_{n+1}  \\
    X^5_{n+3} & X^6_{n+2} & X^7_{n+2}  \\
    \vdots & \vdots & \vdots  \\
    X^{N-3}_{n+l} & X^{N-2}_{n+l} & X^{N-1}_{n+l-1} \\
    X^{N-2}_{n+l+1} & X^{N-1}_{n+l} & X^{N}_{n+l}X^{N+1}_{n+l} \\
\end{pmatrix}.
\] 
The analogue for (\ref{evenDtyperelations}) is
\begin{equation*}
\frac{X^{N+1}_{n+1}+X^{N+1}_{n-1}}{X^N_n}=\frac{X^{N+1}_{n+1}X^N_{n+1}+X^{N-2}_n}{X^{N-1}_{n}}=\frac{X^{N+1}_{n-1}X^N_{n-1}+X^{N-2}_n}{X^{N-1}_{n-1}},
\end{equation*} 
which we use to add the first two black rows to (\ref{oddDfullmatrix}). The following rows are again constructed in the same way we started:
\begin{equation}\label{oddDfullmatrix}
\begin{pmatrix}
\color{grey}X_n^4 & \color{grey}X_n^3 & \color{grey}X_{n-1}^1X_{n-1}^2 \\
\color{grey}X^1_{n+1} & \color{grey}X_n^2 & \color{grey}X^1_{n-1} \\
\color{grey}X^1_{n+1}X^2_{n+1} & \color{grey}X^3_{n+1} & \color{grey}X^4_{n}  \\
\color{grey}X^3_{n+2} & \color{grey}X^4_{n+1} & \color{grey}X^5_{n+1}  \\
\color{grey}X^4_{n+2} & \color{grey}X^5_{n+2} & \color{grey}X^6_{n+1}  \\
\color{grey}X^5_{n+3} & \color{grey}X^6_{n+2} & \color{grey}X^7_{n+2}  \\
\vdots & \vdots & \vdots  \\
\color{grey}X^{N-3}_{n+l} & \color{grey}X^{N-2}_{n+l} & \color{grey}X^{N-1}_{n+l-1} \\
\color{grey}X^{N-2}_{n+l+1} & \color{grey}X^{N-1}_{n+l} & \color{grey}X^{N}_{n+l}X^{N+1}_{n+l} \\
X^{N+1}_{n+l+2} & X^N_{n+l+1} & X^{N+1}_{n+l} \\
X^{N+1}_{n+l+2}X^N_{n+l+2} & X_{n+l+1}^{N-1} & X_{n+l+1}^{N-2} \\
X^{N-1}_{n+l+2} & X^{N-2}_{n+l+2} & X^{N-3}_{n+l+1} \\
\vdots & \vdots & \vdots \\
X^5_{n+2l+1} & X^4_{n+2l} & X^3_{n+2l} \\
X^4_{n+2l+1} & X^3_{n+2l+1} & X^2_{n+2l}X^1_{n+2l} \\
X^1_{n+2l+2} & X^2_{n+2l+1} & X^1_{n+2l}
\end{pmatrix}.
\end{equation} 
From the kernel for this matrix we see that 
\[
J_n:=\frac{X^1_{n+1}+X^1_{n-1}}{X^2_n}
\]
has period $2l+1=N-2$, the same as in the even $N$ case.
\begin{theorem} 
For odd $N$ the linear relation for the extending vertex variables is 
\begin{equation} 
X^k_{n+4N-8}-\mathcal{K}X^k_{n+2N-4}+X^k_n=0.
\end{equation}
\begin{proof} 
As above we have 
\[
\Psi_n:=
\begin{pmatrix}
X^1_{n} & X^2_{n+1}   \\
X^1_{n+N-2} & X^2_{n+N-1}   
\end{pmatrix},
\quad
\tilde{L}_n:=
\begin{pmatrix}
0 & -1   \\
1 & J_{n+1}   
\end{pmatrix}
\quad
\text{ and }
\quad
\Psi_n\tilde{L}_n=
\begin{pmatrix}
X^2_{n+1} & X^1_{n+2}   \\
X^2_{n+N-1} & X^1_{n+N}   
\end{pmatrix}.
\]
Once again we define $L_n:=\tilde{L}_n\tilde{L}_{n+1}$ so $\Psi_nL_n=\Psi_{n+2}$. Applying Theorem \ref{CayleyHamilton} with $q=2$, $p=N-2$ and $m=2N-4$ yields $\Psi_{n+4N-8}-\mathcal{K}\Psi_{n+2N-4}+\Psi_n=0$, from which we can extract the linear relation. The arguments that $\trace(M_n)$ is invariant and $|M_n|=1$ are identical to the ones from the $N$ even case.
\end{proof}
\end{theorem}
The following result is Equation \ref{recurrencesforextendingverticesDodd} in Theorem \ref{similartoAtypetheorem}.
\begin{corollary}
The $X^1_n$ variables satisfy the A type recurrence
\[
X^1_{n+N-1}X^1_n=X^1_{n+N-2}X^1_{n+1}\lambda^2_{n+1}+\gamma^1_n,
\]
where both $\gamma^1_n$ and $\lambda_n$ are period $2$. 
\begin{proof}
By taking determinants in $\Psi_nL_n=\Psi_{n+2}$ we see that
\[
X^1_nX^2_{n+N-1}-X^2_{n+1}X^1_{n+N-2}
\]
is period $2$.
By Lemma \ref{Dtypeperiod2} we can replace $X^2_{n+N-1}$ and $X^2_{n+1}$ to give that
\begin{equation}\label{21stjuly}
X^1_nX^1_{n+N-1}\lambda_n-X^1_{n+1}X^1_{n+N-2}\lambda_{n+1}
\end{equation}
is period $2$, where we have again defined $\lambda_n:=\left(\frac{X^2_0}{X^1_0}\right)^{(-1)^n}$, which is also period $2$. We define expression (\ref{21stjuly}) as $\gamma^1_n/\lambda_{n+1}$ and rearrange to give the corollary.
\end{proof}
\end{corollary}
\subsection{The $\tilde{E}_6$ quiver}
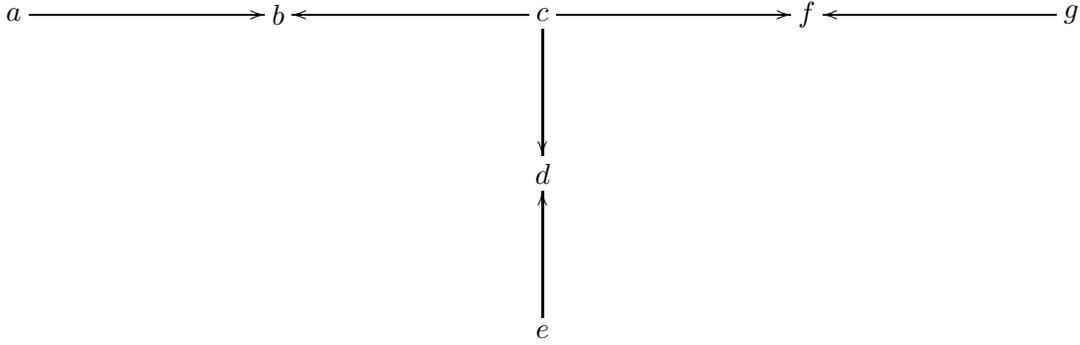
\begin{figure}[h]
\begin{xy} 0;<1pt,0pt>:<0pt,-1pt>:: 
			(0,0) *+{a} ="a",
			(100,0) *+{b} ="b",
			(200,60) *+{d} ="d",
			(200,0) *+{c} ="c",
			(200,120) *+{e} ="e",
			(300,0) *+{f} ="f",
			(400,0) *+{g} ="g",
			"a", {\ar"b"},
			"c", {\ar"d"},
			"c", {\ar"f"},
			"c", {\ar"b"},
			"g", {\ar"f"},
			"e", {\ar"d"},
			
\end{xy}\caption{The $\tilde{E}_6$ quiver.}\label{EuclideanE6}
\end{figure}
The $\tilde{E}_6$ diagram is given in Figure \ref{EuclideanE6} with a source-sink orientation. Since we have a fixed number of vertices we can label each using different letters, making calculations clearer. In our previous notation we would have had, for example, $a_n=X^1_n$. The sequence of mutations, $\mu$, gives the following recurrence relations. 
\begin{align}
a_{n+1}a_{n}&=1+b_n, &
b_{n+1}b_{n}&=1+a_{n+1}c_{n+1}, &
c_{n+1}c_n=1+b_{n}d_{n}f_n, \nonumber\\
d_{n+1}d_n&=1+c_{n+1}e_{n+1}, &
e_{n+1}e_{n}&=1+d_{n}, \nonumber\\
f_{n+1}f_{n}&=1+c_{n+1}g_{n+1}, &
g_{n+1}g_n&=1+f_{n}. \label{mutationrelations}
\end{align}
The diagram has $S_3$ symmetry generated by (for example) the reflections $(bf)(ag)$ and $(bd)(ae)$. It will be useful to name some of the group elements.
\[
\sigma_1:=(bf)(ag), \qquad \sigma_2:=(bd)(ae), \qquad \sigma_3:=(bdf)(aeg)
\]
\subsubsection{Periodic Quantities}
Using the mutation equations, (\ref{mutationrelations}), we form the following matrix, which has determinant $0$ by (\ref{dodgson}). We add labels to the right of the matrix so we may refer to each row:
\begin{equation*}
M:=
\begin{blockarray}{cccc}
\begin{block}{(ccc)c}
1 & a_{n-1} & b_{n-1} & -1 \\
a_n & 1 & 0 & 0\\
b_n & a_{n+1} & 1 & 1 \\
\end{block}
\end{blockarray}
\end{equation*}
After scaling we take a kernel vector $(1,-a_n,1)^T$ since, as before, we can see that the third entry of this vector is equal to the first by (\ref{mutationrelations}). We can add rows to $M$ satisfying the conditions of Theorem \ref{kernelpreservation}:
\begin{equation*}
M:=
\begin{blockarray}{cccc}
\begin{block}{(ccc)c}
\star & d_{n-3}/f_{n-2} & e_{n-2} & -4 \\
b_{n-3} & c_{n-2} & d_{n-2}f_{n-2} & -3 \\
a_{n-2} & b_{n-2} & c_{n-1} & -2 \\
1 & a_{n-1} & b_{n-1} & -1 \\
a_n & 1 & 0 & 0\\
b_n & a_{n+1} & 1 & 1 \\
c_{n+1} & b_{n+1} & a_{n+2} & 2 \\
d_{n+1}f_{n+1} & c_{n+2} & b_{n+2} & 3 \\
e_{n+2} & d_{n+2}/f_{n+1} & \star & 4 \\
\end{block}
\end{blockarray}
\end{equation*}
until we reach rows $\pm 4$ where, in order to use the $d$ mutation relation, we have to divide in our middle entry. The starred entries of $M$ can be filled, preserving the kernel, by insisting that the right $2 \times 2$ minor of rows $3$ and $4$ have determinant $1$, ditto for the left $2 \times 2$ minor of rows $-4$ and $-3$.
\begin{lemma}\label{firstE6lemma}
Setting the rightmost entry of row $4$ to $\frac{g_{n+4}}{f_{n+1}}$ and the leftmost entry of row $-4$ to $\frac{g_{n-4}}{f_{n-2}}$ will preserve the kernel.
\begin{proof}
We just prove the statement for row $4$. In order for the $2 \times 2$ minor to have determinant $1$ we need to set the blank entry to
\[
\frac{1+b_{n+2}d_{n+2}/f_{n+1}}{c_{n+2}}=\frac{f_{n+1}+b_{n+2}d_{n+2}}{c_{n+2}f_{n+1}},
\]
but we have
\[
\frac{f_{n+1}+b_{n+2}d_{n+2}}{c_{n+2}}=\frac{f_{n+1}+\frac{c_{n+3}c_{n+2}-1}{f_{n+2}}}{c_{n+2}}=\frac{f_{n+1}f_{n+2}-1+c_{n+3}c_{n+2}}{c_{n+2}f_{n+2}}
\]
\[
=\frac{c_{n+2}g_{n+2}+c_{n+3}c_{n+2}}{c_{n+2}f_{n+2}}=\frac{g_{n+2}+c_{n+3}}{f_{n+2}}.
\]
Now from the kernel of row $-2$ we have $a_{n-2}-a_nb_{n-2}+c_{n-1}=0$, to which we apply the permutation $\sigma_1$ and shift to see $\frac{g_{n+2}+c_{n+3}}{f_{n+2}}=g_{n+4}$. 
\end{proof}
\end{lemma} 
Using this Lemma we can fill in rows $\pm 4$. We then use vertex $e$ mutation to partially add rows $\pm 5$:
\begin{equation*}
M:=
\begin{blockarray}{cccc}
\begin{block}{(ccc)c}
\star & e_{n-3} & f_{n-2} & -5 \\
g_{n-4}/f_{n-2} & d_{n-3}/f_{n-2} & e_{n-2} & -4 \\
b_{n-3} & c_{n-2} & d_{n-2}f_{n-2} & -3 \\
\vdots & \vdots & \vdots \\
d_{n+1}f_{n+1} & c_{n+2} & b_{n+2} & 3 \\
e_{n+2} & d_{n+2}/f_{n+1} & g_{n+4}/f_{n+1} & 4 \\
f_{n+1} & e_{n+3} & \star & 5 \\
\end{block}
\end{blockarray}
\end{equation*}
Now again we require that both
\[
\begin{vmatrix}
\star & e_{n-3} \\
g_{n-4}/f_{n-2} & d_{n-3}/f_{n-2}
\end{vmatrix}=1,
\qquad
\begin{vmatrix}
d_{n+2}/f_{n+1} & g_{n+4}/f_{n+1} \\
e_{n+3} & \star
\end{vmatrix}=1.
\]
We simply define quantities that satisfy these relations 
\[
J_n:=\frac{f_{n+1}+e_{n+3}g_{n+4}}{d_{n+2}}, \qquad \tilde{J}_n:=\frac{f_{n-2}+e_{n-3}g_{n-4}}{d_{n-3}}
\]
and fill the blank entries with them. The reason for naming these $J_n$ and $\tilde{J}_n$ will be made apparent in the following lemma. We will need the equations from the kernel for rows $\pm 5$:
\begin{equation}\label{Jkernel}
f_{n+1}-a_ne_{n+3}+J_n=0, \qquad \tilde{J}_n-a_ne_{n-3}+f_{n-2}=0.
\end{equation}
\begin{lemma} We have the following expressions for $J_n$ and $\tilde{J}_n$ in terms of the extending variables:
\begin{equation}\label{Js}
J_n=\frac{a_n+g_{n+4}}{e_{n+2}}, \qquad \tilde{J}_n=\frac{a_n+g_{n-4}}{e_{n-2}}
\end{equation}
and each of these has period $3$.
\begin{proof}
From the kernel for row $4$ and vertex $e$ mutation we have 
\[
f_{n+1}=\frac{a_nd_{n+2}-g_{n+4}}{e_{n+2}}=\frac{a_n(e_{n+3}e_{n+2}-1)-g_{n+4}}{e_{n+2}}=a_ne_{n+3}-\frac{a_n+g_{n+4}}{e_{n+2}}.
\]
Comparing this with the first equation of (\ref{Jkernel}) gives $J_n=\frac{a_n+g_{n+4}}{e_{n+2}}$. The corresponding $\tilde{J}_n$ result is proved similarly. Comparing the expressions in (\ref{Js}) one sees that $\sigma_1(J_n)=\tilde{J}_{n+4}$ and applying $\sigma_2$ to (\ref{Jkernel}) we have $\sigma_2(J_{n})=\tilde{J}_{n+3}$. Composing these gives $\sigma_2\sigma_1(J_n)=J_{n+1}$.
The observation that $\sigma_2\sigma_1$ has order $3$ gives the periodicity of $J_n$ and $\tilde{J}_n$. 
\end{proof}
\end{lemma}
\subsubsection{Constant Coefficient Linear Relations}
\begin{theorem}
The recurrences for the variables attached at the extending vertices satisfy the constant coefficient linear relations
\begin{equation}\label{linearrelation}
x_{n+12}-\mathcal{K}x_{n+6}+x_{n}=0,
\end{equation}
i.e. $x\in \{a,e,g\}$ and $\mathcal{K}$ is invariant.
\begin{proof}
Defining
\[ 
\Psi_n:=
\begin{pmatrix}
e_{n+5} & a_{n+3} \\
e_{n+2} & a_{n}
\end{pmatrix},
\qquad
\tilde{L}_n:=
\begin{pmatrix}
J_n & 1 \\
-1 & 0
\end{pmatrix}
\]
we have $\Psi_n\tilde{L}_n=\sigma_3(\Psi_{n+2})$ which gives $\Psi_{n+6}=\Psi_n\tilde{L}_n\tilde{L}_{n+2}\tilde{L}_{n+4}$. We define $L_n:=\tilde{L}_n\tilde{L}_{n+2}\tilde{L}_{n+4}$ and apply Theorem \ref{CayleyHamilton} with $p=3$ and $q=6$. This gives, since $|L_n|=1$,
\[
\Psi_{n+12}-tr(L_n)\Psi_{n+6}-\Psi_n=0.
\]
Note that here $L_n=M_n$. Any entry of this matrix equation will give (\ref{linearrelation}), since $\mathcal{K}:=\trace(L_n)$ is invariant as $\trace(L_{n+1})=\trace(\tilde{L}_{n+1}\tilde{L}_n\tilde{L}_{n+2})=\trace(\tilde{L}_n\tilde{L}_{n+1}\tilde{L}_{n+2})=\trace(L_n)$.
\end{proof}
\end{theorem}
In this case we have slightly more information. We have a period $2$ quantity $K_n$ which we use to find a linear relation for the $b$, $d$ and $f$ vertices. We also prove a relation between $K_n$ and $\mathcal{K}$.
\begin{lemma} The following expressions, each in terms of variables at only one vertex, are fixed by all permutations:
\begin{equation}\label{lemma}
\frac{a_{n-3}+a_{n+3}}{a_n}=\frac{e_{n-3}+e_{n+3}}{e_n}=\frac{g_{n-3}+g_{n+3}}{g_n}.
\end{equation}
\begin{proof}
From (\ref{Jkernel}) we have $J_n-\tilde{J}_{n+3}=a_ne_{n+3}-a_{n+3}e_n$,
but the left hand side is period $3$ so $a_ne_{n+3}-a_{n+3}e_n=a_{n+3}e_{n+6}-a_{n+6}e_{n+3}$ which we rearrange for the result.
\end{proof}
\end{lemma}
The following lemma is equivalent to the $\tilde{E}_6$ part of Theorem \ref{similartoAtypetheorem}.
\begin{lemma}
The quantity
\begin{equation}\label{period2}
K_n:=\frac{a_{n-3}+a_{n+3}}{a_n}
\end{equation}
has period $2$. 
\begin{proof}
Since $J_n=\frac{a_n+g_{n+4}}{e_{n+2}}$ is period $3$ we have
\[
a_ne_{n+5}+e_{n+5}g_{n+4}=a_{n+3}e_{n+2}+e_{n+2}g_{n+7}.
\]
Applying $(\varphi^*)^3+1$ to this equation we get
\[
e_{n+5}(a_n+a_{n+6}+g_{n+4}+g_{n+10})=a_{n+3}(e_{n+2}+e_{n+8})+g_{n+7}(e_{n+2}+e_{n+8}),
\]
which we divide by $a_{n+3}e_{n+5}g_{n+7}$, noting that by (\ref{lemma}) $K_n$ is fixed by permutations,
\[
\frac{K_{n+3}}{g_{n+7}}+\frac{K_{n+7}}{a_{n+3}}=\frac{K_{n+5}}{g_{n+7}}+\frac{K_{n+5}}{a_{n+3}}.
\]
Assuming we may divide by $K_{n+5}-K_{n+3}$ this is equivalent to
\[
\frac{g_{n+7}(K_{n+7}-K_{n+5})}{K_{n+5}-K_{n+3}}=a_{n+3},
\]
but here the left hand side is fixed under $\sigma_2$, so we have $a_{n+3}=e_{n+3}$. This is a contradiction since $a_{n+3}$ and $e_{n+3}$ are independent variables, since setting $n=-3$ means that $a_{n+3}$ and $e-{n+3}$ would appear in a set of arbitrary initial variables. Hence the division by $K_{n+5}-K_{n+3}$ was invalid, giving $K_{n+5}-K_{n+3}=0$.
\end{proof}
\end{lemma}
\begin{lemma} The period $2$ quantity $K_n$ and the invariant $\mathcal{K}$ are related via $\mathcal{K}=K_0K_{1}-2$.
\begin{proof}
From (\ref{period2}) we have
\[
\begin{pmatrix}
a_{n+5} & a_{n+2} \\
a_{n+3} & a_n
\end{pmatrix}
\begin{pmatrix}
K_1 & 1 \\
-1 & 0 
\end{pmatrix}
=
\begin{pmatrix}
a_{n+8} & a_{n+5} \\
a_{n+6} & a_{n+4}
\end{pmatrix}.
\]
Considering the product 
$\begin{pmatrix}
K_1 & 1 \\
-1 & 0 
\end{pmatrix}
\begin{pmatrix}
K_0 & 1 \\
-1 & 0 
\end{pmatrix}$
we apply an argument similar to Theorem \ref{CayleyHamilton} and compare with (\ref{linearrelation}) to arrive at the result.
\end{proof}
\end{lemma}
We can use $K_n$ to find a constant coefficient linear relation for the vertices adjacent to the extending ones. 
\begin{theorem}
We have the linear relation $x_n-(\mathcal{K}+1)x_{n+3}+(\mathcal{K}+1)x_{n+6}-x_{n+9}=0$, for $x\in \{b,d,f\}$.
\begin{proof}
From (\ref{mutationrelations}) We have
\[
\frac{b_n-b_{n+9}}{b_{n+3}-b_{n+6}}=\frac{a_na_{n-1}-a_{n+9}a_{n+8}}{a_{n+3}a_{n+2}-a_{n+6}a_{n+5}}.
\]
Using (\ref{period2}) to replace everything in the numerator on the right hand side, i.e.
\[
a_{n-1}=K_{n+1}a_{n+2}-a_{n+5}, \qquad a_n=K_na_{n+3}-a_{n+6},
\]
\[
a_{n+8}=K_na_{n+5}-a_{n+2}, \qquad a_{n+9}=K_{n+1}a_{n+6}-a_{n+3},
\]
we can rearrange to give
\[
\frac{b_n-b_{n+9}}{b_{n+3}-b_{n+6}}=K_nK_{n+1}-1=\mathcal{K}+1.
\]
This proves the theorem for $x=b$. The other two cases follow by symmetry.
\end{proof}
\end{theorem}
\subsection{The $\tilde{E}_7$ quiver}
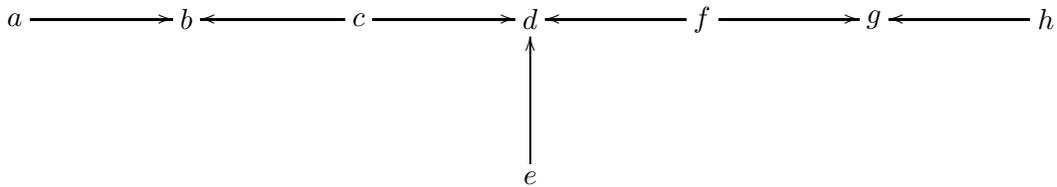
\begin{figure}[h]
\begin{xy} 0;<1pt,0pt>:<0pt,-1pt>::
			(0,0) *+{a} ="a",
			(65,0) *+{b} ="b",
			(130,0) *+{c} ="c",
			(195,60) *+{e} ="e",
			(195,0) *+{d} ="d",
			(260,0) *+{f} ="f",
			(325,0) *+{g} ="g",
			(390,0) *+{h} ="h",
			"a", {\ar"b"},
			"c", {\ar"b"},
			"e", {\ar"d"},
			"f", {\ar"d"},
			"c", {\ar"d"},
			"f", {\ar"g"},
			"h", {\ar"g"},
			
\end{xy}\caption{The $\tilde{E}_7$ quiver.} \label{EuclideanE7}
\end{figure}
The mutation relations for the $\tilde{E}_7$ quiver, Figure \ref{EuclideanE7}, are
\begin{align}
a_{n+1}a_n&=1+b_n, &b_{n+1}b_n&=1+a_{n+1}c_{n+1},  &c_{n+1}c_n&=1+b_nd_n, \nonumber \\
d_{n+1}d_n&=1+c_{n+1}e_{n+1}f_{n+1},  &e_{n+1}e_n&=1+d_n,  &f_{n+1}f_n&=1+d_ng_n, \nonumber \\
g_{n+1}g_n&=1+f_{n+1}h_{n+1},  &h_{n+1}h_n&=1+g_n.
\end{align}
This diagram has $S_2$ symmetry, the generator of which we denote $\sigma:=(ah)(bg)(cf)$.
\subsubsection{Periodic Quantities.} 
We begin with the matrix centred at vertex $a$ with kernel vector $(1,-a_n,1)^T$:
\begin{equation*}
\begin{blockarray}{cccc}
\begin{block}{(ccc)c}
\star & f_{n-3}/e_{n-2} & g_{n-3} & -5 \\
c_{n-3} & d_{n-3} & e_{n-2}f_{n-2} & -4 \\
b_{n-3} & c_{n-2} & d_{n-2} & -3\\
a_{n-2} & b_{n-2} & c_{n-1} & -2 \\
1 & a_{n-1} & b_{n-1} & -1 \\
a_n & 1 & 0 & 0\\
b_n & a_{n+1} & 1 & 1 \\
c_{n+1} & b_{n+1} & a_{n+2} & 2 \\
d_{n+1} & c_{n+2} & b_{n+2} & 3\\
e_{n+2}f_{n+2} & d_{n+2} & c_{n+3} & 4 \\
g_{n+2} & f_{n+3}/e_{n+2} & \star & 5\\
\end{block}
\end{blockarray}
\end{equation*}
The construction of which is routine until we reach the stars in rows $\pm 5$.
\begin{lemma}
Setting the lower right and upper left entries to be
\begin{equation*}
\frac{e_{n+2}+c_{n+3}f_{n+3}}{d_{n+2}e_{n+2}}=\frac{e_{n+4}}{e_{n+2}},
\qquad
\frac{e_{n-2}+f_{n-3}c_{n-3}}{d_{n-3}e_{n-2}}=\frac{e_{n-4}}{e_{n-2}}.
\end{equation*}
respectively preserves the kernel. 
\begin{proof}
The proof of these equalities is the same argument as that used for Lemma \ref{firstE6lemma}.
\end{proof}
\end{lemma}
We can now replace the stars, adding to the rows marked $\pm5$ below. The rows labelled $\pm6$ are to be filled using the mutation relation for $g$, giving us another two starred entries:
\begin{equation*}
\begin{blockarray}{cccc}
\begin{block}{(ccc)c}
\star & g_{n-4} & e_{n-2}h_{n-3} & -6 \\
e_{n-4}/e_{n-2} & f_{n-3}/e_{n-2} & g_{n-3} & -5 \\
\vdots & \vdots & \vdots \\
g_{n+2} & f_{n+3}/e_{n+2} & e_{n+4}/e_{n+2} & 5 \\
e_{n+2}h_{n+3} & g_{n+3} & \star & 6 \\
\end{block}
\end{blockarray}
\end{equation*}
which we can fill with the following result.
\begin{lemma} Due to the relations
\[
\frac{e_{n+2}+e_{n+4}g_{n+3}}{f_{n+3}e_{n+2}}=\frac{a_{n+6}}{e_{n+2}}, \qquad 
\frac{e_{n-2}+e_{n-4}g_{n-4}}{e_{n-2}f_{n-3}}=\frac{a_{n-6}}{e_{n-2}}
\]
we can replace the stars with $a_{n-6}$ and $a_{n+6}$ respectively.
\begin{proof}
From the kernel for row $-5$ we have $e_{n-4}-f_{n-3}a_n+g_{n-3}e_{n-2}=0$, which we shift upwards by $6$ and substitute into our expression to get the first equality. The second is proved similarly.
\end{proof}
\end{lemma}
We now use the mutation relation for $h$ to add most of rows $\pm 7$:
\begin{equation*}
	\begin{blockarray}{cccc}
		\begin{block}{(ccc)c}
			 & h_{n-4}/e_{n-2} & 1 & -7 \\
			a_{n-6} & g_{n-4} & e_{n-2}h_{n-3} & -6  \\
			\vdots & \vdots & \vdots \\
			e_{n+2}h_{n+3} & g_{n+3} & a_{n+6} & 6 \\
			1 & h_{n+4}/e_{n+2} &  & 7 \\
		\end{block}
	\end{blockarray}
\end{equation*}
This is the last row we need to add. We define the quantities 
\[
J_n:=\frac{e_{n+2}+a_{n+6}h_{n+4}}{g_{n+3}}, \qquad \tilde{J}_n:=\frac{e_{n-2}+a_{n-6}h_{n-4}}{g_{n-4}}
\]
that sit in the blank corners, such that the kernel is persevered. As we shall see, these $J_n$ and $\tilde{J}_n$ are also periodic.
\begin{theorem}
$J_n$ and $\tilde{J}_n$ satisfy
\begin{equation}\label{firstJrelation}
J_n=\tilde{J}_{n+6}
\end{equation}
\begin{proof}
Multiplying the numerator and denominator of $J_n$ by $g_{n+2}$ we get
\begin{equation}\label{1}
J_n=\frac{g_{n+2}e_{n+2}+g_{n+2}a_{n+6}h_{n+4}}{g_{n+2}g_{n+3}}.
\end{equation}
From the kernel for row $-6$ we have 
\[
a_ng_{n-4}=a_{n-6}+e_{n-2}h_{n-3},
\]
which we shift upwards by $6$ then multiply by $h_{n+4}$ to give
\begin{equation}\label{2}
a_{n+6}g_{n+2}h_{n+4}=a_{n}h_{n+4}+e_{n+4}h_{n+3}h_{n+4}=a_nh_{n+4}+e_{n+4}(1+g_{n+3}),
\end{equation}
where the second equality uses $h$ mutation. We also have, from the kernel for row $5$,
\begin{equation}\label{3}
e_{n+2}g_{n+2}=a_nf_{n+3}-e_{n+4}.
\end{equation}
Substituting (\ref{2}) and (\ref{3}) into (\ref{1}) gives
\begin{equation}\label{4}
J_n=\frac{a_n(f_{n+3}+h_{n+4})+e_{n+4}g_{n+3}}{g_{n+2}g_{n+3}}.
\end{equation}
From a shift of row $2$ we have $c_{n+3}-a_{n+2}b_{n+3}+a_{n+4}=0$ to which we can apply $\sigma$ to get $f_{n+3}-h_{n+2}g_{n+3}+h_{n+4}=0$. This is used to replace the bracketed term in (\ref{4}), giving 
\[
J_n=\frac{a_nh_{n+2}+e_{n+4}}{g_{n+2}}
\]
which is $\tilde{J}_{n+6}$.
\end{proof}
\end{theorem}
\begin{lemma}\label{lemma1}
The permutation $\sigma$ acts on $J_n$ and $\tilde{J}_n$ by $\sigma(J_n)=\tilde{J}_{n+4}=J_{n+2}$, therefore $J_n$ and $\tilde{J}_n$ are period $4$.
\begin{proof}
From rows $\pm 7$ we have
\[
J_n=a_nh_{n+4}-e_{n+2}, \qquad \tilde{J}_n=a_nh_{n-4}-e_{n-2},
\]
which gives the first equality. Using this we can apply $\sigma$ to (\ref{firstJrelation}) to get $J_{n+2}=\tilde{J}_{n+4}$. 
\end{proof}
\end{lemma}
\begin{lemma} We have the following expression, fixed under $\sigma$, in terms of the variables at only one extending vertex:
\begin{equation}\label{period3quantity}
K_n:=\frac{a_{n+8}+a_n}{a_{n+4}}=\frac{h_{n+8}+h_n}{h_{n+4}}.
\end{equation}
\begin{proof}
Since $J_n=a_nh_{n+4}-e_{n+2}$ and $J_{n+2}=\sigma(J_n)=a_{n+4}h_n-e_{n+2}$ we have
\[
J_n-J_{n+2}=a_nh_{n+4}-a_{n+4}h_n.
\]
The left hand side of this is period $4$ so $a_nh_{n+4}-a_{n+4}h_n$ is also period $4$, which is equivalent to (\ref{period3quantity}).
\end{proof}
\end{lemma}
We have the following more useful expression for $J_n$ which we will use to calculate linear relations.
\begin{theorem} The quantity $J_n$ may be be written
\begin{equation}\label{secondJrelation}
J_n=\frac{a_{n+6}+a_n}{h_{n+3}}=\frac{h_{n+8}+h_{n+2}}{a_{n+5}}.
\end{equation}
\begin{proof}
From row $6$ we see
\begin{equation*}
\begin{split}
0=e_{n+2}h_{n+3}-a_ng_{n+3}+a_{n+6}=e_{n+2}h_{n+3}-a_nh_{n+4}h_{n+3}+a_n+a_{n+6}\\
=-h_{n+3}J_n+a_n+a_{n+6}.
\end{split}
\end{equation*}
We apply $\sigma$ to this equation, noting that by Lemma \ref{lemma1} $\sigma(J_n)=J_{n+2}$, to get the result.
\end{proof}
\end{theorem}
\begin{proposition}
The quantity $K_n$, defined in (\ref{period3quantity}), is period $3$, hence the  $a$ variables satisfy the following linear relation ${a_{n+8}-K_na_{n+4}+a_n=0}$ and the $\tilde{E}_7$ part of Theorem \ref{similartoAtypetheorem}. We also have that the determinant
\[
\tilde{K}_n:=
\begin{vmatrix}
a_n & h_{n+3} \\
a_{n+4} & h_{n+7}
\end{vmatrix}
\]
is period $2$ with $\sigma(\tilde{K}_n)=\tilde{K}_{n+1}$.
\begin{proof}
The proof of this statement is not very illuminating, so we omit it here. It will appear in the author's thesis.
\end{proof}
\end{proposition}
Finally we conjecture an alternate expression for $\tilde{K}_n$.
\begin{conjecture}
The period $2$ quantity $\tilde{K}_n$ can be expressed as
\[
\tilde{K}_n=\frac{a_{n+12}+a_n}{h_{n+6}}.
\]
\end{conjecture}
\subsubsection{Constant Coefficient Linear Relations}
\begin{theorem}
We have the constant coefficient linear relation $x_{n+24}-\mathcal{K}x_{n+12}+x_n=0$, where $x \in \{a,h\}$.
\begin{proof}
We define the matrices
\[
\Psi_n:=
\begin{pmatrix}
a_{n+5} & h_{n+2} \\
h_{n+3} & a_n
\end{pmatrix},
\qquad
\tilde{L}_n:=
\begin{pmatrix}
J_n & 1 \\
-1 & 0 
\end{pmatrix},
\]
then, due to (\ref{secondJrelation}) we have,
\[
\Psi_n\tilde{L}_{n+2}=\sigma(\Psi_{n+3}), \qquad \Psi_n\tilde{L}_{n+2}\tilde{L}_{n+3}=\Psi_{n+6}.
\]
Now we can set $L_n:=\tilde{L}_{n+2}\tilde{L}_{n+3}$ and use Theorem \ref{CayleyHamilton} with $M_n=L_nL_{n+6}$. Again we call $\mathcal{K}:=tr(M_n)$, which, by the same argument used for $\tilde{E}_6$, is invariant.
\end{proof}
\end{theorem}
\subsection{The $\tilde{E}_8$ quiver}
\begin{figure}[h]
\begin{xy} 0;<1pt,0pt>:<0pt,-1pt>::
			(-110,0) *+{a} ="a",
			(-55,0) *+{b} ="b",
			(0,0) *+{c} ="c",
			(55,0) *+{d} ="d",
			(110,0) *+{e} ="e",
			(165,60) *+{g} ="g",
			(165,0) *+{f} ="f",
			(220,0) *+{h} ="h",
			(275,0) *+{i} ="i",
			"a", {\ar"b"},
			"c", {\ar"b"},
			"c", {\ar"d"},
			"e", {\ar"d"},
			"g", {\ar"f"},
			"h", {\ar"f"},
			"e", {\ar"f"},
			"h", {\ar"i"},
\end{xy}\caption{The $\tilde{E}_8$ quiver.} \label{EuclideanE8} 
\end{figure}
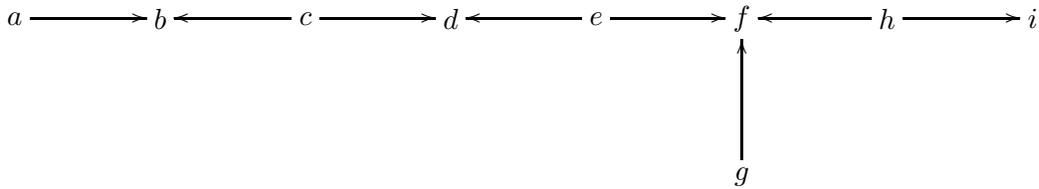
The $\tilde{E}_8$ is given in Figure \ref{EuclideanE8}. The recurrence relations are
\begin{align}
a_{n+1}a_n&=1+b_n,  &b_{n+1}b_n&=1+a_{n+1}c_{n+1},  &c_{n+1}c_n&=1+b_nd_n, \nonumber \\
d_{n+1}d_n&=1+c_{n+1}e_{n+1},  &e_{n+1}e_n&=1+d_nf_n,  &f_{n+1}f_n&=1+e_{n+1}g_{n+1}h_{n+1}, \nonumber \\
g_{n+1}g_n&=1+f_{n},  &h_{n+1}h_n&=1+f_ni_n,  &i_{n+1}i_n&=1+h_{n+1}.
\end{align}
\subsubsection{Periodic Quantities.}
In this subsection we omit some proofs that are similar to previous calculations and will appear in the author's thesis. We begin with the matrix centred at vertex $a$, denoted $M_a$, with kernel vector $(1,-a_n,1)^T$:
\begin{equation*}
M_a:=
\begin{blockarray}{cccc}
\begin{block}{(ccc)c}
\star & h_{n-4}/g_{n-3} & i_{n-4} & -7 \\
e_{n-4} & f_{n-4} & g_{n-3}h_{n-3} & -6 \\
d_{n-4} & e_{n-3} & f_{n-3} & -5 \\
c_{n-3} & d_{n-3} & e_{n-2} & -4 \\
b_{n-3} & c_{n-2} & d_{n-2} & -3\\
a_{n-2} & b_{n-2} & c_{n-1} & -2 \\
1 & a_{n-1} & b_{n-1} & -1 \\
a_n & 1 & 0 & 0\\
b_n & a_{n+1} & 1 & 1 \\
c_{n+1} & b_{n+1} & a_{n+2} & 2 \\
d_{n+1} & c_{n+2} & b_{n+2} & 3\\
e_{n+2} & d_{n+2} & c_{n+3} & 4 \\
f_{n+2} & e_{n+3} & d_{n+3} & 5\\
g_{n+3}h_{n+3} & f_{n+3} & e_{n+4} & 6 \\
i_{n+3} & h_{n+4}/g_{n+3} & \star & 7\\
\end{block}
\end{blockarray}
\end{equation*}
\begin{lemma}\label{glemma} The stars in rows $\pm 7$ are found via the relations
\begin{equation*}
\frac{g_{n+3}+e_{n+4}h_{n+4}}{f_{n+3}g_{n+3}}=\frac{g_{n+5}}{g_{n+3}},
\qquad
\frac{g_{n-3}+e_{n-4}h_{n-4}}{g_{n-3}f_{n-4}}=\frac{g_{n-5}}{g_{n-3}}.
\end{equation*}
The proof of this is analogous to Lemma \ref{firstE6lemma}.
\end{lemma}
We can now fill in  rows $\pm7$ and add some of rows $\pm 8$ via $i$ mutation:
\begin{equation*}
M_a=
\begin{blockarray}{cccc}
\begin{block}{(ccc)c}
\star & i_{n-5} & g_{n-3} & -8 \\
g_{n-5}/g_{n-3} & h_{n-4}/g_{n-3} & i_{n-4} & -7 \\
e_{n-4} & f_{n-4} & g_{n-3}h_{n-3} & -6 \\
\vdots & \vdots & \vdots\\
g_{n+3}h_{n+3} & f_{n+3} & e_{n+4} & 6 \\
i_{n+3} & h_{n+4}/g_{n+3} & g_{n+5}/g_{n+3} & 7\\
g_{n+3} & i_{n+4} & \star & 8 \\
\end{block}
\end{blockarray}
\end{equation*}
\begin{lemma}\label{alemma}
The starred entries can be replaced by $a_{n\pm8}$ respectively, preserving the kernel, because we have
\[
\frac{g_{n+3}+g_{n+5}i_{n+4}}{h_{n+4}}=a_{n+8}, \qquad 
\frac{g_{n-3}+g_{n-5}i_{n-5}}{h_{n-4}}=a_{n-8}.
\]
\begin{proof}
From row $-7$ we have $g_{n-5}-a_nh_{n-4}+g_{n-3}i_{n-4}=0$, which we shift up $8$ and substitute into our expression to get the first result.
\end{proof}
\end{lemma}
We've extended $M_a$ as much as necessary. In this case we'll also need to construct a matrix centred at vertex $i$:
\begin{equation*}
M_i:=
\begin{blockarray}{cccc}
\begin{block}{(ccc)c}
\tilde{J}_n & b_{n-4} & a_{n-3}g_{n-1} & -7 \\
a_{n-8}/g_{n-1} & c_{n-3}/g_{n-1} & b_{n-3} & -6 \\
i_{n-5} & d_{n-3} & c_{n-2}g_{n-1} & -5 \\
g_{n-3}/g_{n-1} & e_{n-2}/g_{n-1} & d_{n-2} & -4 \\
h_{n-2} & f_{n-2} & e_{n-1}g_{n-1} & -3 \\
i_{n-2} & h_{n-1} & f_{n-1} & -2 \\
1 & i_{n-1} & h_n & -1 \\
i_n & 1 & 0 & 0 \\
h_{n+1} & i_{n+1} & 1 & 1 \\
f_{n+1} & h_{n+2} & i_{n+2} & 2 \\
e_{n+2}g_{n+2} & f_{n+2} & h_{n+3} & 3 \\
d_{n+2} & e_{n+3}/g_{n+2} & g_{n+4}/g_{n+2} & 4 \\
c_{n+3}g_{n+2} & d_{n+3} & i_{n+5} & 5 \\
b_{n+3} & c_{n+4}/g_{n+2} & a_{n+9}/g_{n+2} & 6 \\
a_{n+4}g_{n+2} & b_{n+4} & J_n & 7 \\
\end{block}
\end{blockarray}
\end{equation*}
The entries of rows $\pm 4$ and $\pm 5$ were simplified by results similar to Lemmas \ref{glemma} and \ref{alemma} respectively. For rows $\pm 6$ we have use the following equalities:
\[
\frac{g_{n+2}+c_{n+4}i_{n+5}}{d_{n+3}g_{n+2}}=\frac{a_{n+9}}{g_{n+2}}, \qquad \frac{g_{n-1}+c_{n-3}i_{n-5}}{d_{n-3}g_{n-1}}=\frac{a_{n-8}}{g_{n-1}}.
\]
We've also defined 
\[
J_n:=\frac{g_{n+2}+a_{n+9}b_{n+4}}{c_{n+4}}, \qquad \tilde{J}_n:=\frac{g_{n-1}+a_{n-8}b_{n-4}}{c_{n-3}}
\]
The following lemma proves one of the $\tilde{E}_8$ relations given in Theorem \ref{similartoAtypetheorem}.
\begin{lemma}\label{E8lemma}
The expressions $J_n$ and $\tilde{J}_n$ satisfy	$J_n=\tilde{J}_n=a_{n-3}a_{n+4}-i_n$ and $J_n=J_{n+5}$. We also have the following expression which we will use to apply Theorem \ref{CayleyHamilton}:
\begin{equation}\label{period5}
J_{n+3}=\frac{a_{n+12}+a_n}{a_{n+6}}.
\end{equation}
\begin{proof}
We prove only the last relation. We have
\begin{equation*}
\begin{split}
a_{n+12}=c_{n+7}i_{n+3}-b_{n+3}g_{n+2}=a_nb_{n+6}-a_{n+6}i_{n+3}=a_n(-1+a_{n+6}a_{n+7})-a_{n+6}i_{n+3}=\\
-a_n+a_{n+6}(a_na_{n+7}-i_{n+3})=-a_n+a_{n+6}J_{n+3}.
\end{split}
\end{equation*}
The first equality comes from row $6$ of $M_i$ and the second from rows $-2$ and $-8$ of $M_a$. The third uses $a$ mutation and the fifth the first result of this lemma.
\end{proof}
\end{lemma}
To prove the existence of the period $3$ quantity we will use the following matrix, also beginning at vertex $a$, but with alternate rows $7$ and $8$. 
\[
\tilde{M}_a:=
\begin{blockarray}{cccc}
\begin{block}{(ccc)c}
a_n & 1 & 0 & 0 \\
b_n & a_{n+1} & 1 & 1 \\
c_{n+1} & b_1 & a_{n+2} & 2 \\
d_{n+1} & c_{n+2} & b_{n+2} & 3\\
e_{n+2} & d_{n+2} & c_{n+3} & 4\\
f_{n+2} & e_{n+3} & d_{n+3} & 5 \\
g_{n+3}h_{n+3} & f_{n+3} & e_{n+4} & 6 \\
1 & g_{n+4}/h_{n+3} & i_{n+5}/h_{n+3} & 7\\
h_{n+3} & g_{n+4} & i_{n+5} & 8 \\
\end{block}
\end{blockarray}
\]
\begin{lemma}\label{trickylemma}
We have
\[
a_na_{n+13}-i_{n+6}=J_{n+3}J_{n+4}-1.
\]
\begin{proof}
From (\ref{period5}) we have the first equality of the following:
\begin{equation}\label{startoflemma}
a_{n}a_{n+13}=J_{n+4}a_{n}a_{n+7}-a_{n}a_{n+1}=J_{n+4}(J_{n+3}+i_{n+3})-1-b_n. 
\end{equation}
The second equality uses $J_{n+3}=a_na_{n+7}-i_{n+3}$ and $a$ mutation. Now from row $-7$ of $M_i$ we have $i_{n+3}J_{n+4}-i_{n+3}i_{n+4}b_{n}+i_{n+3}a_{n+1}g_{n+3}=0$. We replace the $i_{n+3}i_{n+4}$ term to get
\begin{equation*}
i_{n+3}J_{n+4}-b_{n}=b_{n}h_{n+4}-i_{n+3}a_{n+1}g_{n+3}.
\end{equation*}
The right-hand side of (\ref{startoflemma}) now becomes 
\begin{equation}\label{newRHS}
J_{n+3}J_{n+4}-1+b_nh_{n+4}-i_{n+3}a_{n+1}g_{n+3}.
\end{equation}
Row $7$ of $M_a$ yields 
\[
i_{n+3}a_{n+1}g_{n+3}=a_{n}a_{n+1}h_{n+4}-a_{n+1}g_{n+5}=h_{n+4}+b_nh_{n+4}-a_{n+1}g_{n+5},
\]
so (\ref{newRHS}) is equal to
\begin{equation}\label{newnewRHS}
J_{n+3}J_{n+4}-1+a_{n+1}g_{n+5}-h_{n+4}.
\end{equation}
Finally we have $h_{n+4}-a_{n+1}g_{n+5}+i_{n+6}$ from row $8$ of $\tilde{M}_a$ so (\ref{newnewRHS}) becomes $J_{n+3}J_{n+4}-1-i_{n+6}$ which completes the proof.
\end{proof}
\end{lemma}
The following theorem proves the other $\tilde{E}_8$ relation given in Theorem \ref{similartoAtypetheorem}.
\begin{theorem}
The quantity 
\[
K_n:=\frac{a_{n+20}+a_n}{a_{n+10}}
\]
is period $3$. 
\begin{proof}
By Lemma \ref{trickylemma} the expression $a_na_{n+13}-i_{n+6}$ is period $5$. In particular 
\[
a_na_{n+13}-i_{n+6}=a_{n+10}a_{n+23}-i_{n+16}.
\]
We also have $J_{n+1}=a_{n+3}a_{n+10}-i_{n+6}=a_{n+13}a_{n+20}-i_{n+16}$ so
\[
a_na_{n+13}-a_{n+3}a_{n+10}=a_{n+10}a_{n+23}-a_{n+13}a_{n+20},
\]
which can be factored to give the theorem.
\end{proof}
\end{theorem}
\subsubsection{Constant Coefficient Linear Relations}
\begin{theorem}
The constant coefficient relation for the $a$ variables is $a_{n+60}-\mathcal{K}a_{n+30}+a_n=0$, where $\mathcal{K}$ is invariant.
\begin{proof}
Defining 
\[
\Psi_n:=
\begin{pmatrix}
a_{n+11} & a_{n+5} \\
a_{n+6} & a_n
\end{pmatrix},
\qquad
L_{n-3}:=
\begin{pmatrix}
J_n & 1 \\
-1 & 0
\end{pmatrix} 
\]
gives $\Psi_nL_{n}=\Psi_{n+6}$ so again we apply Theorem \ref{CayleyHamilton} with $p=5$, $q=6$ and ${M_{n}:=L_{n}L_{n+6}\ldots L_{n+24}}$ to arrive at the above result.
\end{proof}
\end{theorem}
Finally, we state that we were unable to prove the following conjecture to complete the table in Figure \ref{periodictable}, but we give an ansatz.
\begin{conjecture}\label{period2conjecture}
The following expression is period $2$:
\[
\tilde{K}_n:=\frac{a_{n+30}+a_n}{a_{n+15}}.
\]
\end{conjecture}
\section{Integrability for $\tilde{D}$ and $\tilde{E}$ quivers}\label{integrabilitysection}
In this section we construct reduced cluster maps, (\ref{reducedclustermap}), and find log-canonical Poisson structures for the reduced variables, (\ref{logcanonicalbracket}). By examining the Poisson subalgebras generated by the $J_n$, and in some cases $K_n$, we find enough commuting first integrals, defined in terms of these periodic quantities, to prove the integrability of the reduced systems. We do this for $\tilde{D}_N$ with $N$ odd and for each $\tilde{E}$ type system. The proof for $\tilde{D}_N$ with $N$ even, however, was too complicated so the problem remains open.
\subsection{Integrability for $\tilde{D}_N$}
The $B$ matrix for the $\tilde{D}$ type diagrams, oriented as above, is
\setcounter{MaxMatrixCols}{20}
\[
B:=\begin{pmatrix}
    0 & 0 & 1 & 0 & 0 & 0 &  &  &  &  &  &  \\
    0 & 0 & 1 & 0 & 0 & 0 &  &  &  &  & &  \\
   -1 & -1 & 0 & -1 & 0 & 0 &  &  &  &  & &  \\
     0 & 0 & 1 & 0 & 1 & 0 &  &  &  &  & & \\
      0 & 0 & 0 & -1 & 0 & -1 & \ddots &  &  & & &    \\
     &  &  & & & \ddots & \ddots & &  &  &  &   \\
    &  &  &  & &  &  & \pm 1 & 0 & \pm 1 & 0 &0 \\
     &  & &  &  & &  & 0 & \mp 1 & 0 & \mp 1 & \mp 1 \\
     &  & &  &  & &  & 0 & 0 & \pm 1 & 0 & 0  \\
     &  & &  &  & &  & 0 & 0 & \pm 1 & 0 & 0 
\end{pmatrix}.
\] 
Where the signs depend on the parity of $N$. The first two rows show that $B$ is singular, so we first project to a lower dimensional space. Since the new coordinates are given in terms of the image of $B$, which is different for different parities of $N$, we'll deal with the two cases separately.

Since the Poisson bracket on the reduced variables commutes with the shift operator $\varphi^*$, calculating the brackets between variables with the same subscripts is straightforward. When calculating the brackets between $J$ variables, however, different subscripts will appear. 
\subsection{The $N$ odd case} 
Let $\{e_i\}$ be standard basis vectors. Here the kernel of $B$ is spanned by $e_1-e_2$ and $e_N-e_{N+1}$  and the image is spanned by 
\[
e_1+e_2, \quad e_3, \quad e_4, \quad \ldots \quad, e_{N-1}, \quad e_N+e_{N+1}.
\]
We take reduced variables $p_n:=X_n^1X_n^2$, $q_n:=X^N_nX^{N+1}_n$ and leave $X^3_n, \ldots X^{N-1}_n$ fixed. This gives a reduced cluster map
\begin{equation}\label{Doddreducedmap}
\varphi:
\begin{pmatrix}
p_n \\
X^3_n \\
\vdots \\
X^{N-1}_n \\
q_n
\end{pmatrix}
\mapsto
\begin{pmatrix}
p_{n+1} \\
X^3_{n+1} \\
\vdots \\
X^{N-1}_{n+1} \\
q_{n+1}
\end{pmatrix},
\end{equation} 
where the $X^i_{n+1}$ are defined as in (\ref{Drelations1}) and (\ref{Drelations3}) for $i=4,\ldots,N-2$. The remaining relations are
\[
p_{n+1}=\frac{1}{p_n}\left(1+X^3_{n+1}\right)^2, \qquad X^3_{n+1}=\frac{1}{X^3_n}\left(1+p_nX^4_n\right),
\]
\[
q_{n+1}=\frac{1}{q_n}\left(1+X^{N-1}_n\right)^2, \qquad X^{N-1}_{n+1}=\frac{1}{X_n^{N-1}}\left(1+q_{n+1}X^{N-2}_{n+1}\right).
\]
Our goal is to prove the integrability of the map (\ref{Doddreducedmap}). As mentioned at the end of Subsection \ref{dynamicalsystemsubsection} the kernel and image of $B$ are used to define $A$:
\[
A:=\begin{pmatrix}
e_1+e_2 \\
e_3\\
\vdots \\
e_{N-1}\\
e_N+e_{N+1}\\
e_1-e_2 \\
e_N-e_{N+1}
\end{pmatrix},
\]
which gives $\hat{B}$ and then $C=(c_{ij})$, the Poisson matrix for the reduced variables:
\[
c_{ij}=\begin{cases}
(-1)^{(j-i+1)/2} & j \text{ even},\: i \text{ odd}\\
0 & \mathrm{otherwise.}
\end{cases}
\]
This gives the Poisson bracket, as in (\ref{logcanonicalbracket}), between $y$ variables, where
\[
y^1_n:=p_n, \quad y^2_n:=X^3_n, \quad y^3_n:=X^4_n, \quad \ldots \quad y^{N-2}_n:=X^{N-1}_n, \quad y^{N-1}_n:=q_n.
\] 
For our calculations, however, we'll stick to using $X$, $p$ and $q$. Since each $J_n$ is preserved by scaling (\ref{lambdascaling}) for each kernel vector of $B$ we can express them in terms of the reduced variables, hence we can calculate brackets between the $J_n$.  
\begin{theorem}
The Poisson structure for the $J$ variables, defined in Lemma \ref{periodiclemma}, is given by
\begin{equation}\label{oddNDpoisson}
\{J_0,J_{-k}\}=(-1)^kJ_0J_{-k}+\delta_{1,k}-\delta_{N-3,k},
\end{equation}
for $k=1,...,N-3$, where $\delta_{i,j}$ is the Kronecker delta. We can apply the shift $\varphi^*$, a Poisson algebra homomorphism, to this relation to obtain the brackets between the remaining $J$ variables.
\begin{proof}
We use (\ref{oddDfullmatrix}) to express $J_n$ in plenty of different ways, i.e. for each row
\[
J_n=\frac{1\textsuperscript{st} \textrm{ entry} + 3\textsuperscript{rd} \textrm{ entry} }{2\textsuperscript{nd} \textrm{ entry}}.
\] 
Shifting these expressions by appropriate amounts allows us to write each $J_n$ in terms of reduced variables with subscript $0$ or $1$, simplifying calculations. We have
\[
J_0=\frac{p_1+X^4_0}{X^3_1}, \qquad J_{\frac{-N+3}{2}}=\frac{X^{N-2}_1+q_0}{X_0^{N-1}}
\]
and for $k=1,\ldots ,\frac{N-5}{2}$
\[
J_{-k}=\frac{X_1^{2k+2}+X_0^{2k+4}}{X_1^{2k+3}}.
\]
We first calculate 
\[
\{J_0,J_{-k}\}=\left\{\frac{p_1+X^4_0}{X^3_1},\frac{X_1^{2k+2}+X_0^{2k+4}}{X_1^{2k+3}}\right\}=\frac{1}{X^3_1X_1^{2k+3}}\{p_1+X^4_0,X_1^{2k+2}+X_0^{2k+4}\}
\]
\[
-\frac{J_{-k}}{X^3_1X^{2k+3}_1}\{p_1+X^4_0,X_1^{2k+3}\}-\frac{J_0}{X_1^{2k+3}X^3_1}\{X^3_1,X_1^{2k+2}+X_0^{2k+4}\}.
\]
The issue now is calculating the brackets between elements with unequal subscripts. For example
\[
\{p_1,X^{2k+4}_0\}=\left\{p_1,\frac{1+X^{2k+3}_1X^{2k+5}_1}{X^{2k+4}_1}\right\}=\frac{1}{X^{2k+4}_1}\{p_1,X^{2k+3}_1X^{2k+5}_1\}=
\]
\[
\frac{X^{2k+3}_1}{X^{2k+4}_1}\{p_1,X^{2k+5}_1\}+\frac{X^{2k+5}_1}{X^{2k+4}_1}\{p_1,X^{2k+3}_1\}=0.
\]
After some work we have
\[
\{J_0,J_{-k}\}=\delta_{k,1}+(-1)^kJ_0J_{-k},
\]
while the bracket
\[
\{J_0,J_{\frac{N-3}{2}}\}=(-1)^{(N-3)/2}J_0J_{\frac{N-3}{2}}
\]
is found similarly.
\end{proof}
\end{theorem}
\begin{theorem}
The reduced cluster map (\ref{Doddreducedmap}) is Liouville integrable. 
\begin{proof}
Equation (\ref{oddNDpoisson}) can be written as the sum of the two Poisson brackets
\[
\{J_0,J_{-k}\}_0:=\delta_{1,k}-\delta_{N-3,k}, \qquad \{J_0,J_{-k}\}_2:=(-1)^kJ_0J_{-k}.
\] 
This bi-Hamiltonian structure \cite{magri} allows for the construction of $(N-1)/2$ commuting first integrals. The Poisson subalgebra (\ref{oddNDpoisson}) is the same, up to scaling, as that found for the $\tilde{A}_{1,N}$ quivers in \cite{fordy}, where the authors give a more complete proof. See \cite{fordy2011mutation} for further details.
\end{proof}
\end{theorem}
\subsection{The case $N=6$}
For general even $N$ the calculations were too much to bear. Instead we will work with an example with $N=6$. The $B$ matrix and our choice of $A$ are
\[
B=
\begin{pmatrix}
0 & 0 & 1 & 0 & 0 & 0 & 0 \\
0 & 0 & 1 & 0 & 0 & 0 & 0 \\
-1 & -1 & 0 & -1 & 0 & 0 & 0 \\
0 & 0 & 1 & 0 & 1 & 0 & 0 \\
0 & 0 & 0 & -1 & 0 & -1 & -1 \\
0 & 0 & 0 & 0 & 1 & 0 & 0 \\
0 & 0 & 0 & 0 & 1 & 0 & 0
\end{pmatrix},
\qquad A=
\begin{pmatrix}
0 & 0 & 1 & 0 & 0 & 0 & 0 \\
1 & 1 & 0 & 1 & 0 & 0 & 0 \\
0 & 0 & 0 & 1 & 0 & 1 & 1 \\
0 & 0 & 0 & 0 & 1 & 0 & 0 \\
1 & -1 & 0 & 0 & 0 & 0 & 0 \\
0 & 0 & 0 & 0 & 0 & 1 & -1 \\
1 & 0 & 0 & -1 & 0 & 0 & 1
\end{pmatrix}.
\]
This gives the Poisson matrix
\[
C=
\begin{pmatrix}
0 & 1 & 0 & 0 \\
-1 & 0 & 0 & 0 \\
0 & 0 & 0 & -1 \\
0 & 0 & 1 & 0
\end{pmatrix}.
\]
Our reduced variables are 
\[
y^1_n:=p_n:=X_n^1X_n^2X_n^4, \quad y^2_n:=X_n^3, \quad  y^3_n:=X_n^5, \quad y^4_n:=q_n:=X_n^4X_n^6X_n^7,
\]
giving a reduced cluster map
\begin{equation}\label{D=6ddreducedmap}
\varphi:
\begin{pmatrix}
p_n \\
X^3_n \\
X^{5}_n \\
q_n
\end{pmatrix}
\mapsto
\begin{pmatrix}
p_{n+1} \\
X^3_{n+1} \\
X^{5}_{n+1} \\
q_{n+1}
\end{pmatrix},
\end{equation} 
where
\begin{equation*}
p_{n+1}=\frac{1}{p_n}\left(1+X^3_{n+1}\right)^2\left(1+X^3_{n+1}X^5_{n+1}\right), \qquad X^3_{n+1}=\frac{1}{X^3_n}\left(1+p_n\right)
\end{equation*}
\[
X^5_{n+1}=\frac{1}{X^5_n}=\left(1+q_n\right), \qquad q_{n+1}=\frac{1}{q_n}\left(1+X^3_{n+1}X^5_{n+1}\right)\left(1+X^5_{n+1}\right)^2.
\]
Now we can express the $J$ variables as, for example
\begin{equation}\label{N=6expressionsforJ}
J_n=\frac{X^1_{n+1}X^2_{n+1}+X^4_n}{X^3_{n+1}}, \qquad J_{n-1}=\frac{X^4_{n+1}+X^6_{n}X^7_n}{X^5_{n+1}}.
\end{equation}
The scaling (\ref{lambdascaling}), for the kernel vector $(1,0,0,-1,0,0,1)^T$, acts as
\[
X^1_{n+1}\mapsto \lambda X^1_{n+1}, \quad X^4_{n+1}\mapsto \lambda^{-1}X^4_{n+1},\quad X^7_{n+1}\mapsto \lambda X^7_{n+1}
\]
and fixes the remaining $X^i_{n+1}$. From the expression $X^4_{n+1}X^4_n=1+X^3_{n+1}X^5_{n+1}$ and that of $J_n$ in (\ref{N=6expressionsforJ}) we see that
\[
X^4_n\mapsto \lambda X^4_n, \qquad J_n\mapsto \lambda J_n.
\]
Since $J_n$ is not preserved by this scaling it cannot be written in terms of the reduced variables. However, one sees from the second expression in (\ref{N=6expressionsforJ}), that $J_{n-1}$ transforms to $\lambda^{-1}J_{n-1}$. We instead define $J'_n:=J_nJ_{n-1}$ and look at the Poisson subalgebra generated by $\{J'_{n-i}\}_{i=0}^3$. One can check that these generators are preserved by the scaling associated with the other two kernel vectors of $B$, so can be written in terms of the reduced variables. Indeed, for example, we have
\[
J'_n=\frac{1}{X^3_{n+1}X^5_{n+1}}\left(p_{n+1}+q_n+1+X^3_{n+1}X^5_{n+1}+\frac{p_{n+1}}{q_{n+1}}(1+X^5_{n+1})^2 \right).
\]
We find similar expressions for the other $J'$ and rewrite so that the same subscript appears throughout. With the aid of a computer we arrive at
\[
\{J'_0,J'_{-1}\}=-J'_0J'_{-1}+J'_0+J'_{-1}, \qquad \{J'_{0},J'_{-2}\}=J'_{-3}-J'_{-1}.
\]
We remark that this Poisson subalgebra structure also appears in the $6$ dimensional reduced system from the $\tilde{E}_7$ quiver. Here however, unlike the $\tilde{E}_7$ case, the $J_n$ generated enough independent, commuting first integrals to prove integrability:
\[
J'_0+J'_1+J'_2+J'_3, \qquad J'_0J'_1J'_2J'_3.
\] 
\subsection{Integrability for $\tilde{E}_6$}
For the $\tilde{E}_6$ quiver, we have the $B$ matrix
\begin{equation*}
B=
\begin{pmatrix}
0 & 1 & 0 & 0 & 0 & 0 & 0 \\
-1 & 0 & -1 & 0 & 0 & 0 & 0 \\
0 & 1 & 0 & 1 & 0 & 1 & 0 \\
0 & 0 & -1 & 0 & -1 & 0 & 0 \\
0 & 0 & 0 & 1 & 0 & 0 & 0 \\
0 & 0 & -1 & 0 & 0 & 0 & -1 \\
0 & 0 & 0 & 0 & 0 & 1 & 0
\end{pmatrix}.
\end{equation*} 
For which we can take image vectors 
\[
e_1+e_3,e_2,e_3+e_5,e_4,e_3+e_6,e_7,
\] 
giving reduced coordinates 
\[
y^1_n:=a_nc_n, \qquad y^2_n:=b_n, \qquad y^3_n:=c_ne_n,
\]
\[
y^4_n:=d_n, \qquad y^5_n:=c_ng_n, \qquad y^6_n:=f_n.
\]
The reduced cluster map
\begin{equation}\label{E=6ddreducedmap}
\varphi:
\begin{pmatrix}
y^1_n \\
y^2_n \\
\vdots \\
y^6_n
\end{pmatrix}
\mapsto
\begin{pmatrix}
y^1_{n+1} \\
y^2_{n+1} \\
\vdots \\
y^6_{n+1}
\end{pmatrix}
\end{equation} 
is given by
\[
y^1_{n+1}=\frac{1}{y^1_n}\left(1+y^2_n\right)\left(1+y^2_ny^4_ny^6_n\right), \qquad y^2_{n+1}=\frac{1}{y^2_n}\left(1+y^1_{n+1}\right), \qquad y^3_{n+1}=\frac{1}{y^3_n}\left(1+y^2_ny^4_ny^6_n\right)\left(1+y^4_n\right)
\]
\[
y^4_{n+1}=\frac{1}{y^4_n}\left(1+y^3_{n+1}\right), \qquad y^5_{n+1}=\frac{1}{y^5_n}\left(1+y^2_ny^4_ny^6_n\right)\left(1+y^6_n\right), \qquad y^6_{n+1}=\frac{1}{y^6_n}\left(1+y^5_{n+1}\right).
\]
The image vectors and a kernel vector $(e_1-e_3+e_5+e_7)$ are used as a basis for the matrix $A$. We may then calculate the Poisson matrix, after scaling,
\[
C=
\begin{pmatrix}
0 & 1 & 0 & 0 & 0 & 0 \\
-1 & 0 & 0 & 0 & 0 & 0 \\
0 & 0 & 0 & 1 & 0 & 0 \\
0 & 0 & -1 & 0 & 0 & 0 \\
0 & 0 & 0 & 0 & 0 & 1 \\
0 & 0 & 0 & 0 & -1 & 0
\end{pmatrix}.
\]
We can explicitly write $J_0$ as 
\[
J_0=\frac{y^1_0y^5_0+y^1_0y^6_0+y^1_0+y^5_0+(1+y^6_0)(y^2_0y^4_0y^6_0+1)}{y^2_0y^5_0y^6_0}
\]
and find similar expressions for each of $J_1,J_2,\tilde{J}_0,\tilde{J}_1,\tilde{J}_2$ from the $S_3$ action. The first row of the Poisson matrix is
\[
\{J_0,J_1\}=J_0J_1-1, \qquad \{J_0,J_2\}=-J_0J_2+1, \qquad \{J_0,\tilde{J_i}\}=0
\]
for $i=1,2,3$. From this we can find each bracket, noting that the tilde is an automorphism. The algebra is the direct sum of the subalgebras generated by the $J_i$ and the $\tilde{J}_i$ and these summands are isomorphic. We take independent commuting first integrals
\[
J_0+J_1+J_2, \qquad J_0J_1J_2, \qquad \tilde{J}_0+\tilde{J}_1+\tilde{J}_2,
\]
proving the integrability of (\ref{E=6ddreducedmap}). 
\subsection{Integrability for $\tilde{E}_7$}\label{E_7IntegrabilitySection}
In this case we have the $B$ matrix
\begin{equation*}
B=
\begin{pmatrix}
0 & 1 & 0 & 0 & 0 & 0 & 0 & 0 \\
-1 & 0 & -1 & 0 & 0 & 0 & 0 & 0 \\
0 & 1 & 0 & 1 & 0 & 0 & 0 & 0 \\
0 & 0 & -1 & 0 & -1 & -1 & 0  & 0\\
0 & 0 & 0 & 1 & 0 & 0 & 0 & 0 \\
0 & 0 & 0 & 1 & 0 & 0 & 1 & 0 \\
0 & 0 & 0 & 0 & 0 & -1 & 0 & -1 \\
0 & 0 & 0 & 0 & 0 & 0 & 1 & 0 
\end{pmatrix}.
\end{equation*} 
Our choice of image vectors are 
\[
e_1+e_3,e_2,e_3+e_5+e_6,e_4,e_6+e_8,e_7
\]
and reduced coordinates
\[
y^1_n:=a_nc_n, \qquad y^2_n:=b_n, \qquad y^3_n:=c_ne_nf_n,
\]
\[
y^4_n:=d_n, \qquad y^5_n:=f_nh_n, \qquad y^6_n:=g_n.
\]
Here the reduced cluster map is given by
\[
y^1_{n+1}=\frac{1}{y^1_n}\left(1+y^2_n\right)\left(1+y^2_ny^4_n\right), \qquad y^2_{n+1}=\frac{1}{y^2_n}\left(1+y^1_{n+1}\right), \qquad y^3_{n+1}=\frac{1}{y^3_n}\left(1+y^2_ny^4_n\right)\left(1+y^4_n\right)\left(1+y^4_ny^6_n\right)
\]
\[
y^4_{n+1}=\frac{1}{y^4_n}\left(1+y^3_{n+1}\right)\, \qquad y^5_{n+1}=\frac{1}{y^5_n}\left(1+y^4_{n}y^6_{n}\right)\left(1+y^6_{n}\right), \qquad y^6_{n+1}=\frac{1}{y^6_n}\left(1+y^5_{n+1}\right).
\]
The matrix $B$ has kernel vectors $e_1-e_3+e_6-e_8$ and $e_1-e_3+e_5$ so we construct $A$ as usual and find
\[
C=
\begin{pmatrix}
0 & 1 & 0 & 0 & 0 & 0 \\
-1 & 0 & 0 & 0 & 0 & 0 \\
0 & 0 & 0 & 1 & 0 & 0 \\
0 & 0 & -1 & 0 & 0 & 0 \\
0 & 0 & 0 & 0 & 0 & 1 \\
0 & 0 & 0 & 0 & -1 & 0
\end{pmatrix}.
\]
The kernel vectors give the scaling symmetries
\[
\lambda_1:
\begin{pmatrix}
a_n \\
b_n \\
c_n \\
d_n \\
e_n \\
f_n \\
g_n \\
h_n
\end{pmatrix}
\mapsto
\begin{pmatrix}
\lambda^{(-1^n)}a_n \\
b_n \\
\lambda^{(-1^{n+1})}c_n \\
d_n \\
e_n \\
\lambda^{(-1^{n})}f_n \\
g_n \\
\lambda^{(-1^{n+1})}h_n
\end{pmatrix},
\qquad 
\lambda_2:
\begin{pmatrix}
a_n \\
b_n \\
c_n \\
d_n \\
e_n \\
f_n \\
g_n \\
h_n
\end{pmatrix}
\mapsto
\begin{pmatrix}
\lambda^{(-1^n)}a_n \\
b_n \\
\lambda^{(-1^{n+1})}c_n \\
d_n \\
\lambda^{(-1^{n})}e_n \\
f_n \\
g_n \\
h_n
\end{pmatrix}.
\]
Our periodic quantities are fixed by $\lambda_1$ but not by $\lambda_2$.  We have
\[
\lambda_2:
\begin{pmatrix}
J_0 \\
J_1 \\
J_2 \\
J_3
\end{pmatrix}
\mapsto
\begin{pmatrix}
\lambda J_0 \\
\lambda^{-1}J_1 \\
\lambda J_2 \\
\lambda^{-1} J_3
\end{pmatrix},
\]
so we can instead take $J'_i:=J_iJ_{i+1}$ for $i=0,1,2$ with each $J'_i$ fixed under $\lambda_1$ and $\lambda_2$. We have, for example,
\[
J'_0=J_0J_1=\frac{1}{y^2_0y^3_0(y^4_0)^2y^5_0y^6_0}(y^1_0y^3_0+y^1_0y^4_0y^6_0+y^3_0+y^1_0+(y^4_0y^6_0+1)(y^2_0y^4_0+1))\times
\]  
\[
(y^4_0(y^4_0((y^6_0)^2+y^6_0)+y^5_0+(y^6_0)^2+2y^6_0+1)+(y^5_0+y^6_0+1)(y^3_0+1).
\]
Now we can calculate
\[
\{J'_0,J'_1\}=J'_0J'_1-J'_0-J'_1, \qquad \{J'_0,J'_2\}=J'_1-J'_3,
\] 
from which the other brackets follow. We have commuting first integrals
\[
J'_0+J'_1+J'_2+J'_3, \qquad  J'_0J'_2+J'_1J'_3.
\]
Similar calculations give that the first integral $K_0+K_1+K_2$ is independent of, and commutes with, the previous two. This proves integrability of the reduced cluster map for $\tilde{E}_7$.
\subsection{Integrability for $\tilde{E}_8$}
Here the $B$ matrix is
\begin{equation*}
B=
\begin{pmatrix}
0 & 1 & 0 & 0 & 0 & 0 & 0 & 0 & 0 \\
-1 & 0 & -1 & 0 & 0 & 0 & 0 & 0 & 0 \\
0 & 1 & 0 & 1 & 0 & 0 & 0 & 0 & 0 \\
0 & 0 & -1 & 0 & -1 & 0 & 0  & 0 & 0\\
0 & 0 & 0 & 1 & 0 & 1 & 0 & 0 & 0 \\
0 & 0 & 0 & 0 & -1 & 0 & -1 & -1 & 0 \\
0 & 0 & 0 & 0 & 0 & 1 & 0 & 0 & 0  \\
0 & 0 & 0 & 0 & 0 & 1 & 0 & 0 & 1 \\
0 & 0 & 0 & 0 & 0 & 0 & 0 & -1 & 0
\end{pmatrix},
\end{equation*} 
with image vectors
\[
e_1+e_3,e_2,e_3+e_5,e_4,e_5+e_7+e_8,e_6,e_8,e_9
\]
and kernel vector $e_1-e_3+e_5-e_7$. These give the matrix
\[
C=
\begin{pmatrix}
0 & 1 & 0 & 0 & 0 & 0 & 0 & 0 \\
-1 & 0 & 0 & 0 & 0 & 0 & 0 & 0 \\
0 & 0 & 0 & 1 & 0 & 0 & 0 & 0\\
0 & 0 & -1 & 0 & 0 & 0 & 0 & 0 \\
0 & 0 & 0 & 0 & 0 & 1 & 0 & 0 \\
0 & 0 & 0 & 0 & -1 & 0 & 0 & 0 \\
0 & 0 & 0 & 0 & 0 & 0 & 0 & 1 \\
0 & 0 & 0 & 0 & 0 & 0 & -1 & 0
\end{pmatrix}.
\]
The reduced coordinates are
\[
y^1_n:=a_nc_n, \quad y^2_n:=b_n, \quad y^3_n:=c_ne_n, \quad y^4_n:=d_n,
\] 
\[
y^5_n:=e_ng_nh_n, \quad y^6_n:=f_n, \quad y^7_n:=h_n, \quad y^8_n:=i_n,
\]
with reduced cluster map $\varphi$, given by
\[
y^1_{n+1}=\frac{1}{y_n^1}\left(1+y^2_n\right)\left(1+y^2_ny^4_n\right), \quad y^2_{n+1}=\frac{1}{y^2_n}\left(1+y^1_{n+1}\right), \quad y^3_{n+1}=\frac{1}{y_n^3}\left(1+y^2_ny^4_n\right)\left(1+y^4_ny^6_n\right),
\]
\[
y^4_{n+1}=\frac{1}{y^4_n}\left(1+y^3_{n+1}\right), \qquad y^5_{n+1}=\frac{1}{y_n^5}\left(1+y^4_ny^6_n\right)\left(1+y^6_n\right)\left(1+y^6_ny^8_n\right),
\]
\[
y^6_{n+1}=\frac{1}{y^6_n}\left(1+y^5_{n+1}\right), \quad y^7_{n+1}=\frac{1}{y^7_n}\left(1+y^6_{n}y^8_{n}\right), \quad y^8_{n+1}=\frac{1}{y^8_n}\left(1+y^7_{n+1}\right).
\]
The single kernel vector generates the scaling
\[
\lambda:
\begin{pmatrix}
a_n \\
b_n \\
c_n \\
d_n \\
e_n \\
f_n \\
g_n \\
h_n \\
i_n
\end{pmatrix}
\mapsto
\begin{pmatrix}
\lambda^{(-1^n)}a_n \\
b_n \\
\lambda^{(-1^{n+1})}c_n \\
d_n \\
\lambda^{(-1^n)}e_n \\
f_n \\
\lambda^{(-1^n)}g_n \\
h_n \\
i_n
\end{pmatrix},
\]
for which the period $5$ quantity
\[
J_n=\frac{a_{n+12}+a_n}{a_{n+6}}
\]
is fixed. Again we compute the brackets
\[
\{J_0,J_1\}=J_0J_1-1, \qquad \{J_0,J_2\}=-J_0J_2,
\] 
and find the three commuting first integrals
\[ 
J_0+J_1+J_2+J_3+J_4, \qquad J_0J_1J_2J_3J_4, \qquad J_0J_1J_2+J_1J_2J_3+J_2J_3J_4+J_3J_4J_0+J_4J_0J_1.
\]
From the period $3$ quantities we construct a final first integral $K_0+K_1+K_2$ which commutes with the other three, proving integrability. 
\bibliographystyle{plain}
\bibliography{AffineDETypes}
\end{document}